\documentclass[3p,sort&compress,times]{elsarticle}
\usepackage[T1]{fontenc}
\usepackage[utf8]{inputenc}
\usepackage[english]{babel}
\usepackage{amsmath} 
\usepackage{pdflscape}
\usepackage{pgfplots}
\pgfplotsset{compat=1.13}
\usepgfplotslibrary{groupplots}
\usepackage{booktabs}
\usepackage{hyperref}
\usepackage{url}
\usepackage{filecontents}
\usepackage{color}
\newcommand{\calT}{\mathcal{T}}
\newcommand{\calV}{\mathcal{V}}
\newcommand{\diff}{\mathrm{d}}
\newcommand{\dx}{\diff x}
\newcommand{\dd}[2][{}]{\frac{\diff#1}{\diff#2}}
\newcommand{\ddt}{\dd[]t}
\newcommand{\pp}[2][{}]{\frac{\partial#1}{\partial#2}}
\newcommand{\sech}{\mathop{\mathrm{sech}}\nolimits}
\newcommand{\eg}{\emph{e.g.}}
\newcommand{\ie}{\emph{i.e.}}

\pgfplotsset{
  compat=1.13,
  every axis plot/.append style={thick},
  /pgf/declare function={
    A(\z,\a,\b) = 
      1-(0.5-(\a*\b*(1-\a-\b)-\a^2*\b^2*(1-2*\b)*(1-2*\a)/2*\z)*\z)*\z;
    degen(\x) = \x/(6*\x-1);
    discrp(\x) = (\x*(4*\x-3)+sqrt(2*\x*(2*\x - 1)^3))/(4*\x*(3*\x-2));
    discrm(\x) = (\x*(4*\x-3)-sqrt(2*\x*(2*\x - 1)^3))/(4*\x*(3*\x-2));
    eff0(\x) = sqrt((((\x-2)*\x+1)*\x-1/6)*\x);
    effp(\x) = ((1-\x)+eff0(\x)/\x)/2;
    effm(\x) = ((1-\x)-eff0(\x)/\x)/2;
    a0(\x) = sqrt((\x*(\x*(\x-2)+1)+1/2)/\x);
    ap(\x) = (\x-1+a0(\x))/2;
    am(\x) = (\x-1-a0(\x))/2;
    b0(\x) = a0(-\x);
    bm(\x) = (\x+1-b0(\x))/2 ;
  },
}

\hypersetup{
  unicode=true,
  breaklinks=true,
  colorlinks=true,
  linkcolor=blue,
  citecolor=blue,
  filecolor=blue,
  urlcolor=MidnightBlue,
  pdftitle={Palindromic 3-stage splitting integrators, a roadmap},
  pdfauthor={Cédric M. Campos, J. M. Sanz-Serna},
  pdfsubject={PACS 02.50.Ng, 02.60.Lj, 02.70.Ns; MSC2010 60J22, 65C40, 65L05, 70F99},
  pdfkeywords={splitting algorithms, Verlet integrator, molecular dynamics, partial differential equations, Hamiltonian Monte Carlo},
  pdfcreator={GNU Emacs 24.5.1},
}

\begin{document}
\title{Palindromic 3-stage splitting integrators, a roadmap\tnoteref{mcode}}
\tnotetext[mcode]{MATLAB implementations of some utilitarian scripts can be found at \url{http://github.com/vitaminace33}.}
\author[yt,uva]{Cédric~M. Campos\corref{cor}}
\ead{cedricmc@yachaytech.edu.ec}
\author[c3]{J.~M. Sanz-Serna}
\ead{jmsanzserna@gmail.com}
\cortext[cor]{Corresponding author}
\address[yt]{Dept. Matemáticas, Universidad Yachay Tech, Hda. San José s/n, 100115 Urcuquí, Ecuador}
\address[c3]{Dept. Matemáticas, Universidad Carlos III de Madrid, Avenida de la Universidad 30, E-28911 Leganés (Madrid), Spain}
\address[uva]{Dept. Matemática Aplicada e IMUVA, Universidad de Valladolid, Paseo de Belén 7, 47011 Valladolid, Spain}

\begin{abstract}
The implementation of multi-stage splitting integrators is essentially the same as the implementation of the familiar Strang/Verlet method. Therefore multi-stage formulas may be easily incorporated into software that now uses the Strang/Verlet integrator. We study in detail the two-parameter family of palindromic, three-stage splitting formulas and identify choices of parameters that may outperform the Strang/Verlet method. One of these choices leads to a method of effective order four suitable to integrate in time some partial differential equations. Other choices may be seen as perturbations of the Strang method that increase efficiency in molecular dynamics simulations and in Hybrid Monte Carlo sampling.
\end{abstract}


\maketitle

\section{Introduction}
We study in detail the two-parameter family of palindromic, three-stage splitting formulas and identify choices of parameters that may outperform the Strang/Verlet method.

The implementation of multi-stage splitting integrators is essentially the same as the implementation of the familiar Strang/Verlet method and, therefore, multi-stage formulas may be easily incorporated into software that now uses the Strang/Verlet integrator. Unfortunately multi-stage integrators do not appear to be as widely used as they deserve; a possible reason for such a lack of popularity is the difficulty associated with choosing the parameters that specify the integrator. For instance, when dealing with the two-parameter family of palindromic,\footnote{The use of the word palindromic to refer to splitting methods like \eqref{eq:meth} that read the same from left to right as from right to left is not universal in the literature. The alternative terms self-adjoint and symmetric are used in the monographs \cite{ssc, hlw} and \cite{bc} respectively. However, as explained in \cite{sirev}, there are reasons to avoid the use of adjoint in this context and of course the word symmetric is overused in science.} three-stage formulas considered in this paper, it is tempting to choose the parameters so as to achieve order four; this is very likely to be an ill-advised decision. Our aim here is to describe in detail the plane of the parameters \((a,b)\) in the family under consideration so as to identify \lq good\rq\ choices of coefficients. These choices have to be tailored to the problem at hand, \ie\ there is no \lq universally best\rq\ choice, and we give specific recommendations in this regard.

Section~\ref{sec:family} describes the family under consideration. Sections~\ref{sec:order} and \ref{sec:harmonic} are devoted to the analysis; the former studies the behavior of the methods as the step size \(h\) approaches 0 and the latter presents finite \(h\) results for the harmonic oscillator model problem. Numerical experiments are reported in Section~\ref{sec:numerics}. We successively consider: (i) the time integration of partial differential equations when the spatial errors are small, (ii) the same issue for large spatial errors, (iii) molecular dynamics simulations, (iv) Hybrid (or Hamiltonian) Monte Carlo sampling. In the final section we make specific recommendations on the choice of parameters, tailored to the problem under consideration.

Splitting integrators have become widely used in recent years. Since they may take advantage of specific features of the system being integrated, they often outperform \lq general purpose\rq\ methods like Runge-Kutta schemes. They are of much interest in connection with geometric integration in general, \ie with the preservation of geometric features of the solution \cite{geometric,hlw}, and in particular with the symplectic integration of Hamiltonian problems \cite{ssc}. Of course they are also relevant in the time-integration of partial differential equations as illustrated in Section~\ref{sec:numerics}. An excellent general reference for splitting algorithms is \cite{bc}, where,  among other things, the reader may find many useful algorithms and a good coverage of the literature.

Throughout the paper much attention is given to the idea of processing, an inexpensive way of enhancing the performance of some integrators, as it appears not to be as well known as it deserves. In particular we show the advantages of processed integrators in the time integration of some partial differential equations. Processing was introduced by J.C. Butcher in 1972 \cite{Bu72}. At that time there was much interest in developing software packages for ordinary differential equations  and the technical difficulties of combining processing with variable step sizes precluded Butcher's ideas from becoming popular. The idea of processing was revitalized later \cite{processing,hessian} in the context of geometric integration where constant step sizes are the rule \cite{ssc}.

\section{The family of palindromic 3-stage splitting methods}
\label{sec:family}
In this section we present the integrators studied in the paper.

\subsection{Splitting algorithms}
Let the system of differential equations of interest be
\begin{equation}\label{eq:AB}
\ddt x = F(x), \qquad F(x) = A(x)+B(x),
\end{equation}
where it is assumed that the \emph{split systems}
\begin{eqnarray}\label{eq:A}
\ddt x = A(x),\\
\label{eq:B}
\ddt x = B(x),
\end{eqnarray}
may both be integrated in closed form. We denote by \(\varphi_t^{(A)}\) the solution flow of \eqref{eq:A}, \ie, for each \(x_0\) and \(t\), \(\varphi_t^{(A)}(x_0)\) is the value at time \(t\) of the solution of \eqref{eq:A} with initial value \(x(0)=x_0\). Similarly, \(\varphi_t^{(B)}\) stands for the solution flow of \eqref{eq:B}. As an illustration, consider the familiar Newton equations for \(x = (q,p)\)
\begin{equation}\label{eq:newton}
\ddt \left[\begin{matrix}q\\p\end{matrix}\right] = \left[\begin{matrix}M^{-1}p\\0\end{matrix}\right]+\left[\begin{matrix}0\\f(q)\end{matrix}\right]
\end{equation}
(\(q\), \(p\), \(f\) are the vectors of coordinates, momenta, and forces respectively and \(M\) is the constant matrix of masses). Here
\[
\varphi_t^{(A)}(q,p) = (q+tM^{-1}p, p ),\qquad \varphi_t^{(B)}(p,q) = (q,p+tF(q)).
\]
These are some times called \lq drift\rq\ and \lq kick\rq\ respectively.

For future reference we note that solution flows satisfy the group property:
\begin{equation}\label{eq:flow}
\varphi_{t_1+t_2}^{(A)} = \varphi_{t_1}^{(A)}\circ \varphi_{t_2}^{(A)},\qquad
\varphi_{t_1+t_2}^{(B)} = \varphi_{t_1}^{(B)}\circ \varphi_{t_2}^{(B)},
\end{equation}
where \(\circ\) denotes the composition of maps, \ie\ \(\varphi_{t_1}^{(A)}\circ \varphi_{t_2}^{(A)}\) is the map that at \(x\) takes the value
\(\varphi_{t_1}^{(A)}\big(\varphi_{t_2}^{(A)}(x)\big)\), etc.

If \(h>0\) is the step size and \(x_n\) the numerical solution at time \(nh\), splitting integrators for \eqref{eq:AB} are of the form \(x_{n+1} =\psi_h(x_n)\), where the mapping \(\psi_h\) is obtained by concatenating flows of the split systems. The best known example is the Strang splitting \cite{strang}
\begin{equation}\label{eq:vv}
\psi_h=\varphi_{(1/2)h}^{(B)}\circ \varphi_{h}^{(A)}\circ\varphi_{(1/2)h}^{(B)}
\end{equation}
with second order accuracy (\ie\ the error in a single step is \(\mathcal{O}(h^3)\) so that, after \(n\) steps, \(x_n\) differs from the true \(x(nh)\) by \(\mathcal{O}(h^2)\) for \(nh\) ranging in a bounded interval). In the particular case of the Newton equations \eqref{eq:newton}, Strang's formula yields the well-known \emph{velocity Verlet} integrator, the method of choice in molecular dynamics \cite{schlick,leim}. It is of course possible to interchange the roles of A and B in \eqref{eq:vv}, \ie\ 
\begin{equation}\label{eq:pv}
\psi_h = \varphi_{(1/2)h}^{(A)}\circ \varphi_{h}^{(B)}\circ\varphi_{(1/2)h}^{(A)};
\end{equation}
for the Newton equations, one then obtains the \emph{position Verlet} integrator \cite{schlick}.

When \( N\) steps of the method \eqref{eq:vv} are taken, the map that advances from \(x_0\) to \(x_N\), \ie\ 
\[ \psi_h^N = \overbrace{
\Big(\varphi_{(1/2)h}^{(B)}\circ \varphi_{h}^{(A)}\circ\varphi_{(1/2)h}^{(B)}\Big)
\circ \cdots\circ
\Big(\varphi_{(1/2)h}^{(B)}\circ \varphi_{h}^{(A)}\circ\varphi_{(1/2)h}^{(B)}\Big)
}^{N \:\textrm{times}}
\]
may be rewritten with the help of \eqref{eq:flow} in the leapfrog form
\[
\psi_h^N = \varphi_{(1/2)h}^{(B)}\circ \overbrace{\Big( \varphi_{h}^{(A)}\circ\varphi_{h}^{(B)}\Big) \circ \cdots\circ
\Big( \varphi_{h}^{(A)}\circ\varphi_{h}^{(B)}\Big)
}^{N-1\: \textrm{times}}\circ\varphi_{h}^{(A)}\circ
\varphi_{(1/2)h}^{(B)};
\]
now the right hand-side only uses \(N+1\) times the flow \(\varphi_t^{(B)}\). Consequently, in the velocity Verlet integrator, the value of the force \(f\) to be used in the second kick of the current step may be reused in the first kick of the next step, so that \(N+1\) force evaluations are sufficient to perform \(N\) integration steps. Since in a single time-step, the use of \eqref{eq:vv} essentially requires one computation of the flow of \eqref{eq:A} and one computation of the flow of \eqref{eq:B}, we say that the scheme has \emph{one stage}.

\subsection{Palindromic 3-stage splittings}
In this article we consider splitting algorithms that in a single time-step use several computations of the flows of \eqref{eq:A} and \eqref{eq:B}. Specifically we study the two-parameter family of 3-stage palindromic splitting formulas
defined by
\begin{equation}\label{eq:meth}
\psi_h =
\varphi_{(1/2-b)h}^{(B)}\circ \varphi_{ah}^{(A)}\circ\varphi_{bh}^{(B)}\circ
\varphi_{(1-2a)h}^{(A)}\circ
\varphi_{bh}^{(B)}\circ\varphi_{ah}^{(A)}\circ \varphi_{(1/2-b)h}^{(B)}.
\end{equation}
The parameters \(a\) and \(b\) are supposed to be real, but the complex case has also been considered \cite{alexander,philippe}. In some applications, the split flow \(\varphi^{(A)}_t\) (or \(\varphi^{(B)}_t\)) only makes sense if \(t\geq 0\). For instance, if the split system A (or B) models diffusion, it cannot be integrated with \(t<0\) as the backward integration of diffusion equations is an ill posed problem. The same difficulty may arise if A (or B) involves Brownian noise.

As explained in the case of Strang's splitting, the last \(B\)-flow of the current time-step may be combined with the first of the next step, so that, per step, \eqref{eq:meth} essentially involves the computation of three \(A\)-flows and three \(B\)-flows. Since one step with a method of the family \eqref{eq:meth} costs three times as much as one of \eqref{eq:pv}, fair comparisons have to run \eqref{eq:meth} with a step length three times longer than the step length one would use for \eqref{eq:pv}. It is obvious that one may consider families of methods with four, five, \dots\ stages \cite{bc}; those methods possess many free parameters and have to be used with very long values of \(h\) if they have to compete with
\eqref{eq:pv} and will not be considered in the present study.

The parametrization used in \eqref{eq:meth} has some degeneracies, in the sense that for some parameter values the number of stages \eqref{eq:meth} is actually lower than three.

\begin{enumerate}
\item For \( a = 0\), the second and sixth mappings in the right hand-side of \eqref{eq:meth} are the identity and we may write
\[\psi_h=
\varphi_{(1/2-b)h}^{(B)}\circ \varphi_{bh}^{(B)}\circ
\varphi_{h}^{(A)}\circ
\varphi_{bh}^{(B)}\circ \varphi_{(1/2-b)h}^{(B)} = \varphi_{(1/2)h}^{(B)}\circ
\varphi_{h}^{(A)}\circ \varphi_{(1/2)h}^{(B)};
\]
thus \(\psi_h\) is the Strang splitting \eqref{eq:vv}, regardless of the value of \(b\).
\item For \( b = 0\), \eqref{eq:meth} reduces similarly to \eqref{eq:vv}, regardless of the value of \(a\).
\item For \( a = 1/2\),
\[
\psi_h = \varphi_{(1/2-b)h}^{(B)}\circ \varphi_{(1/2)h}^{(A)}\circ
\varphi_{2bh}^{(B)}\circ\varphi_{(1/2)h}^{(A)}\circ \varphi_{(1/2-b)h}^{(B)}
\]
and we have a one-parameter family of two-stage integrators.
\item For \( b = 1/2\), we have, symmetrically,
\[
\psi_h = \varphi_{ah}^{(A)}\circ \varphi_{(1/2)h}^{(B)}\circ
\varphi_{h}^{(A)}\circ\varphi_{(1/2)h}^{(B)}\circ \varphi_{ah}^{(A)}.
\]
\end{enumerate}

When both \( a=1/2\) and \(b=1/2\), \eqref{eq:meth} reduces to the one-stage \eqref{eq:pv}. For the choice \(a=b=1/3\), one step of length \(h\) with \eqref{eq:meth} yields exactly the same result as three successive steps of \eqref{eq:vv} each of length \(h/3\). In addition, when \(a = 1/2\) and \(b= 1/4\), one step with \eqref{eq:meth} reproduces two steps of length \(h/2\) with \eqref{eq:vv}. Symmetrically, for \(a = 1/4\) and \(b= 1/2\), \eqref{eq:meth} reproduces two steps of length \(h/2\) with \eqref{eq:pv}.

The square \( 0<a<1/2\), \( 0<b<1/2\) in the \((a,b)\) plane corresponds to methods where all the time-increments used in \eqref{eq:meth} (namely \((1/2-b)h\), \(ah\), \(bh\), \((1/2-a)h\)) are strictly positive; these methods have wider applicability than methods with negative time-increments as we shall discuss later. Methods on the boundary of the square are degenerate.

Multistage splitting integrators may be successfully incorporated into molecular dynamics packages, see \eg\ \cite{elena}, as they can be implemented with the same ease as the Strang/Verlet scheme. (This remark is valid even in constrained cases where Verlet becomes SHAKE or RATTLE.)

In the particular case of the Newton equations \eqref{eq:newton} the integrators in \eqref{eq:meth} coincide with Runge-Kutta-Nystr\"{o}m (RKN) methods \cite{ssc}. Some integrators and results discussed below for the problem \eqref{eq:AB} first appeared in the RKN literature for the Newton equations.

\subsection{Geometric properties}
Splitting integrators qualify as geometric integrators \cite{geometric,hlw} due to some favorable properties:
\begin{itemize}
\item \emph{Preservation of volume.} Assume that, in \eqref{eq:AB}, \(F\), \( A\), \(B\), are divergence-free so that for, each \(t\), the flows \(\varphi_t^{(F)}\), \(\varphi_t^{(A)}\), \(\varphi_t^{(B)}\) preserve volume (\ie\ the determinant of the corresponding Jacobian matrices is 1). Then, for each \(h\), \(\psi_h\) in \eqref{eq:meth} is also a volume-preserving transformation, because the composition of volume-preserving maps is itself volume-preserving.
\item \emph{Symplecticity.} Assume that, in \eqref{eq:AB}, the vector field \(F\) is Hamiltonian with Hamiltonian function \(H_F\), \ie\ the differential system is of the form
\[
\ddt q^i = +\pp[H_F]{p^i}(p,q), \qquad \ddt p^i = -\pp[H_F]{q^i}(p,q)
\]
(\(q^i\) and \(p^i\) are the components of \(x = (q,p)\)). Furthermore assume that the split systems \eqref{eq:A} and \eqref{eq:B} are also Hamiltonian with Hamiltonian functions \(H_A\) and \(H_B\) respectively. Then the flows \(\varphi_t^{(F)}\), \(\varphi_t^{(A)}\), \(\varphi_t^{(B)}\) and the mapping \(\psi_h\) in \eqref{eq:meth} are symplectic transformations. The Newton equations \eqref{eq:newton} provide an example if the forces \(f(q)\) are conservative, \ie\ if \(f(q)=-\nabla V(q)\); then \(H_A\) and \(H_B\) coincide with the kinetic and potential energies \((1/2) p^TM^{-1}p\) and \(V(q)\) respectively.
\item \emph{Reversibility.} Assume that \(F\), \( A\), \(B\) are reversible with respect to a (constant) involution matrix \(S\). This means that \(S^2 = I\) and \(F(Sx) = -SF(x)\), \(A(Sx) = -SA(x)\), \(B(Sx) = -SB(x)\) and implies the reversibility of the flows: \(S\varphi_t^{(F)}(x) = \big(\varphi_t^{(F)}\big)^{-1}(Sx)\), etc. Then the methods in \eqref{eq:meth} are also reversible: \(S\psi_h(x) = \big(\psi_h\big)^{-1}(Sx)\), a property that follows from their palindromic structure. It is well known that the Newton equations \eqref{eq:newton} are reversible with respect to the involution \((q,p)\mapsto (q,-p)\) (time-reversibility of mechanics).
\end{itemize}

\section{Accuracy}
\label{sec:order}
This section is devoted to investigating the accuracy of the integrators of the family \eqref{eq:meth}.

\subsection{Local error}

Although alternative methodologies are available (see \cite{ordercond,AlSa16}), the standard way to analyze the local error of splitting integrators is as follows \cite{ssc,hlw,bc}. Each of the flows being composed is represented as the exponential of a differential operator (Lie derivative) and then the exponentials are combined, by means of the Baker-Campbell-Hausdorff formula, into a single exponential. In this way one finds a \emph{modified} differential equation (parameterized by \(h\))
\begin{equation}\label{eq:modsys}
\ddt x= \tilde F_h(x),
\end{equation}
\ie\ an equation whose \(h\)-flow formally coincides with the map \(\psi_h\) \cite{modified}. Thus, the difference between \(\tilde F_h(x)\) and \(F(x)\) measures the error in a single step of the integrator starting from \(x\). For a method of the form \eqref{eq:meth}, one finds in this way
\begin{equation}\label{eq:Fh}
\tilde F_h = A + B + h^2\alpha\: [A,[A,B]] + h^2 \beta\: [B,[A,B]]+\mathcal{O}(h^4),
\end{equation}
with
\[
\alpha = a^2b-\frac{1}{24},\qquad \beta = -ab^2+ab-\frac{1}{12}.
\]
Here \([\cdot,\cdot]\) represents the commutator or Lie-Jacobi bracket:
\[
[F,G] = \pp[F]{x} G- \pp[G]{x} F
\]
(with an alternative definition, many authors reverse the signs in the right hand-side). If the right hand-side of \eqref{eq:Fh} is expanded further, there are \emph{six} \(h^4\) terms, like \(h^4\gamma [A,[A,[A,[A,B]]]]\),
\(h^4\delta [B,[A,[A,[A,B]]]]\), etc. Note also that, as expected, \eqref{eq:Fh} is \emph{not} symmetric in \(A\) and \(B\). By interchanging the roles of \(A\), \(B\), each splitting algorithm may always be applied in two different forms to a given split system, and these two forms will produce different errors.

From \eqref{eq:Fh}, the order of the method is, at least, two. To achieve order four, we solve the equations \(\alpha =\beta=0\); there is only one solution, often associated with Yoshida's name \cite{yoshida} (but see \cite{candy,forest} for earlier appearances of these coefficients):
\begin{equation}\label{eq:yoshida}
\textrm{Yoshida:}\qquad a = \frac16 \left( 1- \sqrt[3]{2}-\frac{1}{\sqrt[3]{2}} \right) ,\qquad b = 1-2a.
\end{equation}

In the Hamiltonian case, the symplecticity of the integrator implies that the modified system \eqref{eq:modsys} is itself Hamiltonian \cite{fields}; the corresponding Hamiltonian function
\(H_{\tilde F_h}\) satisfies, in view of \eqref{eq:Fh},
\begin{equation}\label{eq:Hh}
H_{\tilde F_h} = H_A + H_B + h^2\alpha\: \{H_A,\{H_A,H_B\}\} + h^2 \beta\: \{H_B,\{H_A,H_B\}\}+\mathcal{O}(h^4),
\end{equation}
where the Poisson bracket of Hamiltonian functions is defined as
\[
\{H, K\} = \sum_{j=1}^d \left( \pp[H]{q_i} \pp[K]{p_i}-\pp[H]{p_i} \pp[K]{q_i}\right).
\]
(Recall that the commutator of two Hamiltonian vector fields it itself a Hamiltonian vector field; with our choice of signs in the commutator, the Hamiltonian of the commutator is given by the Poisson bracket of the
Hamiltonians of the fields being commuted \cite{arnold}.)

\subsection{Local energy error} 
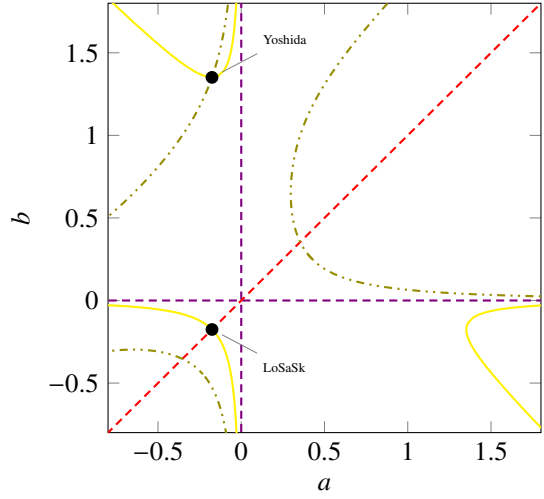
\begin{figure}
\centering
\begin{tikzpicture}
\begin{axis}[
  xmin=-0.8,
  xmax= 1.8,
  ymin=-0.8,
  ymax= 1.8,
  xlabel={\(a\)},
  ylabel={\(b\)},
  unit vector ratio*=1 1 1,
  samples=50,
]

\addplot[violet,densely dashed] coordinates {(-2.0, 0.0) (3.0,0.0)};
\addplot[violet,densely dashed] coordinates {( 0.0,-2.0) (0.0,3.0)};
\addplot[red,densely dashed] coordinates {(-2,-2) (3,3)};

\addplot[yellow,domain=-1.5:-0.18]{effp(x)};
\addplot[yellow,domain=-1.5:-0.18](effp(x),x);
\addplot[yellow,domain=1.36:3]{effp(x)};
\addplot[yellow,domain=-1.5:-0.13](effm(x),x);
\addplot[yellow,domain=1.36:3](effp(x),x);
\addplot[yellow,domain=-1.5:-0.13]{effm(x)};

\addplot[olive,dash dot dot,domain=0.2:3](am(x),x);
\addplot[olive,dash dot dot,domain=0.25:3](ap(x),x);
\addplot[olive,dash dot dot,domain=0.40:3]{bm(x)};
\addplot[olive,dash dot dot,domain=-1.5:-0.23]{bm(x)};
\addplot[olive,dash dot dot,domain=-1.5:-0.44](ap(x),x);

\addplot[mark=*] coordinates {(-0.1756,-0.1756)} node[pin=-30:{\tiny LoSaSk}]{};
\addplot[mark=*] coordinates {(-0.1756, 1.3512)} node[pin=30:{\tiny Yoshida}]{};
\end{axis}
\end{tikzpicture}
\caption{The \((a,b)\) plane. Points on the olive double-dot-dashed curve yield methods with enhanced conservation of energy for quadratic Hamiltonians. Points on the yellow solid curve give methods with effective order four. The intersection of the two curves yields Yoshida's method. The intersection of the yellow solid curve with the red dashed diagonal \(a=b\) yields the LoSaSk method}
\label{fig:order}
\end{figure}

In the Hamiltonian case, the solutions of \eqref{eq:AB} preserve the value of the Hamiltonian function (energy) \(H_F\) and it is of interest in many applications to study to what extent the mapping \(\psi_h\) also preserves that value. Since \(\psi_h\) coincides with the solution of the \(h\)-flow of the modified Hamiltonian, standard results from the Hamiltonian formalism show that 
\[
H_F(\psi_h(x)) =  \exp(h\{\cdot,H_{\tilde F_h}\} H_F(x)= H_F(x)+h\{H_F,H_{\tilde F_h}\}+\cdots
\]
and therefore the energy increment
\[
\Delta= \Delta(x)= H_F(\psi_h(x)) - H_F(x)
\]
may be computed as
\[
\Delta = h\{H_F,H_{\tilde F_h}\}+\cdots
\]
By using \eqref{eq:Hh}, we find, at leading order,
\[ \Delta = h^3\alpha\: \{H_A,\{H_A,\{H_A,H_B\}\}\}+ h^3(\alpha+\beta)\: \{H_A,\{H_B,\{H_A,H_B\}\}\} + h^3\beta\: \{H_B,\{H_B,\{H_A,H_B\}\}\}\}+\cdots \]
This expression shows that for \emph{general} Hamiltonian functions, the energy changes by an p\(\mathcal{O}(h^3)\) amount at each step, except for the choice of \(a\) and \(b\) in Yoshida \eqref{eq:yoshida} where the change is \(\mathcal{O}(h^5)\). However, in particular cases these estimates may be improved. For a separable quadratic Hamiltonian with
\begin{equation}\label{eq:quadratic}
H_A = \frac{1}{2} p^TM^{-1}p,\qquad H_B = \frac{1}{2} p^TK^{-1}p
\end{equation}
(\(K\) is the constant stiffness matrix), the iterated brackets \(\{H_A,\{H_A,\{H_A,H_B\}\}\}\) and \(\{H_B,\{H_B,\{H_A,H_B\}\}\}\}\) vanish and, accordingly, at leading order
\[ \Delta = h^3(\alpha+\beta)\: \{H_A,\{H_B,\{H_A,H_B\}\}\}+\cdots \]
By imposing the condition
\begin{equation}\label{eq:alphaminusbeta}
\alpha=-\beta
\end{equation}
we find a one-parameter family of integrators with enhanced energy conservation \(\Delta=\mathcal{O}(h^5)\) for the separable, quadratic case. These methods correspond to the olive double-dot-dashed curve in Fig.~\ref{fig:order}, with equation
\[
a^2b-ab^2+ab-\frac{1}{8} = 0.
\]

\subsection{Processing. Effective order}
Splitting schemes (and other one-step integrators) may be used in combination with \emph{processing} \cite{processing,hessian}. If \(\psi_h\) is the integrator to be processed, the mapping that advances one step with the processed method is defined as
\begin{equation}\label{eq:psihat}
\widehat \psi_h = \chi_h^{-1}\circ\psi_h\circ\chi_h,
\end{equation}
where the pre-processor \(\chi_h\) is a suitably chosen map and the post-processor \(\chi_h^{-1}\) is the map that inverts the action of \(\chi_h\). If \(N\) steps of the processed method are taken, we may write
\[
\widehat{\psi}_h^N = \overbrace{\Big(\chi_h^{-1}\circ \psi_h\circ\chi_h\Big)\circ \cdots \circ\Big(\chi_h^{-1}\circ \psi_h\circ\chi_h\Big)}^{N\: \textrm{times}}= \chi_h^{-1}\circ \overbrace{\psi_h\circ \cdots \circ\psi_h}^{N\: \textrm{times}}\circ\ \chi_h.
\]
It follows that the pre-processor needs to be applied only at the initial point \(x_0\) and the post-processor at those steps where output is required; as a consequence the work required to implement \(\widehat \psi_h\) is typically not significantly larger than that required to implement the unprocessed \(\psi_h\). The advantage of processing is that the processed integrator may have higher accuracy than \(\psi_h\).

In an alternative description, in processing one generates a numerical trajectory \(X_{n+1} = \psi_h(X_n)\), \(n=0, 1,\dots\), and then performs a change of variables \(X_n=\chi_h(x_n)\); the points \(x_n\) are seen as the approximations to the true \(x(nh)\).

In this paper, methods of the family \eqref{eq:meth}, are processed by taking \(\chi_h\) to be the \(h\) flow of a system of the form
\begin{equation}\label{eq:lambda}
\ddt x = h \lambda\, [A,B]
\end{equation}
(the commutator is evaluated at \(x\)). Note that then the change of variables \(\chi_h\) differs from the identity in \(\mathcal{O}(h^2)\) terms. In the Hamiltonian case, \eqref{eq:lambda} is the Hamiltonian system with Hamiltonian \(h\lambda\{H_A,H_B\}\) and accordingly \(\chi_h\) is a symplectic transformation. It follows from \eqref{eq:psihat} that the processed method is also symplectic. Similarly \(\widehat{\psi}_h\) is volume-preserving for divergence-free \(F\), \(A\), \(B\), because the commutator of divergence-free fields is itself divergence-free. Unfortunately, in the reversible case, it is easily seen that \(\widehat{\psi}_h\) fails to be reversible; implications of this fact will be relevant in Section~\ref{sec:numerics} below.

The modified system of the processed method is found by changing variables in the modified system of \(\psi_h\); in this way we find
\[
\widetilde{\widehat{F_h}}= A + B + h^2\widehat{\alpha}\: [A,[A,B]] + h^2 \widehat{\beta}\: [B,[A,B]] + \mathcal{O}(h^4),
\]
with
\[
\widehat\alpha = \alpha-\lambda = a^2b-\frac{1}{24}-\lambda, \qquad \widehat\beta = \beta-\lambda=-ab^2+ab-\frac{1}{12}-\lambda.
\]
Hence methods with
\begin{equation}\label{eq:alphabeta}
\alpha=\beta\end{equation} may be processed with \(\lambda = \alpha =\beta\), to get \(\widehat\alpha = \widehat\beta = 0\) and then \(\widehat{\psi}_h\) will have order four. One says that in this case \(\psi_h\) has \emph{effective order} four \cite{effective}.

The solid curve in Fig.~\ref{fig:order}, with equation
\[
ab(a+b-1)+\frac{1}{24} = 0
\]
corresponds to integrators with effective order four. The solid and dashed line only intersect at the point defined in \eqref{eq:yoshida} (Yoshida's method). The intersection of the solid curve and the diagonal \(a=b\) corresponds to the method
\begin{equation}\label{eq:LoSaSk}
\textrm{LoSaSk:}\qquad a = \frac16 \left( 1- \sqrt[3]{2}-\frac{1}{\sqrt[3]{2}} \right) ,\qquad b = a,
\end{equation}
that will be discussed later (LoSaSk is an acronym for the names of the authors of \cite{maximal}).

We observe that, for all methods with effective order four, one of the increments \(ah\), \((1/2-a)h\) is negative \emph{and} one of the increments \(bh\), \((1/2-b)h\) is also negative; therefore methods with effective order four are only applicable when \emph{both} split systems may be integrated forward and backward (see \cite{negative}).

It is also possible to process in more sophisticated ways. For instance, we could choose \(\chi_h\) to be the \(h\)-flow of
\begin{equation}\label{eq:sophisticated}
\ddt x = h \lambda [A,B]+\mu [A,[A,B]]+\nu [B,[A,B]],
\end{equation}
rather than of \eqref{eq:lambda}. For methods with effective order four, the parameter \(\lambda\) is determined as before, that is \(\lambda = \alpha = \beta\), in order to annihilate the \(\mathcal{O}(h^2)\) terms of the modified \(\widetilde{\widehat{F_h}}\); the additional parameters \(\mu\) and \(\nu\) are then chosen to minimize the (six) coefficients of the \(\mathcal{O}(h^4)\) terms of \(\widetilde{\widehat{F_h}}\). See \cite{processing,maximal,hessian}, but those references only consider the case of the Newton equations where some of the six iterated commutators with five letters vanish.

\section{The linear model problem}
\label{sec:harmonic}

The (scalar) harmonic oscillator
\[
\ddt q = p,\qquad \ddt p = -q,
\]
provides the simplest example of the Newton equations \eqref{eq:newton} and has been traditionally used as a model problem to discriminate among different integrators. Results obtained for the harmonic oscillator are easily extended to all quadratic problems \eqref{eq:quadratic} via diagonalization (this of course includes many problems resulting from the spatial discretization of Hamiltonian partial differential equations). In addition they often provide some indication on the performance of the method in more complicated, nonlinear cases. Unlike the analysis in the preceding section, that only takes into account the behavior of the integrators as \(h\rightarrow0\), the material below looks at situations with \(h\) fixed.

\subsection{The numerical solution}
The harmonic oscillator may be used to illustrate the role of the conditions \eqref{eq:alphaminusbeta} and \eqref{eq:alphabeta} related to improved conservation of energy and effective order four respectively. The true solution is given by the rotation
\begin{equation}\label{eq:rotation}
\left[ \begin{matrix}q(t)\\p(t)\end{matrix}\right] = M_t\left[ \begin{matrix}q(0)\\p(0)\end{matrix}\right],\qquad
M_t = \left[ \begin{matrix}\phantom{-}\cos t & \sin t\\ -\sin t & \cos t\end{matrix}\right].
\end{equation}
A time-step $(q_{n+1},p_{n+1})= \psi_h(q_n,p_n)$ with a method of the family \eqref{eq:meth} is given by
\begin{equation}\label{eq:harmonicintegrator}
\left[ \begin{matrix}q_{n+1}\\p_{n+1}\end{matrix}\right] = \tilde{M}_h\left[ \begin{matrix}q_n\\p_n\end{matrix}\right],\qquad
\tilde{M}_h=
\left[ \begin{matrix}A_h& B_h\\ C_h & A_h\end{matrix}\right],
\end{equation}
for suitable coefficients $A_h$, $B_h$, $C_h$, that depend on \(h\) and on the parameters \(a\) and \(b\) that specify the method. For instance,
\begin{equation}\label{eq:ah}
A_h(z) = 1-\tfrac12 z+ab(1-a-b)z^2-2a^2b^2(\tfrac12-a)^2(\tfrac12-b)^2z^3, \qquad z = h^2.
\end{equation}
The symplecticity of the integrators implies that \(A_h^2-B_hC_h=1\).

The evolution over $n$ time-steps is then
\begin{equation}\label{eq:harmonicintegrator2}
\left[ \begin{matrix}q_{n}\\p_{n}\end{matrix}\right] = \tilde{M}_h^n\left[ \begin{matrix}q_0\\p_0\end{matrix}\right].
\end{equation}
Assuming that \(|A_h|<1\) (which is equivalent to \(M_h\) possessing a pair of complex conjugate eigenvalues of unit modulus), we introduce $\theta_h$ such that $A_h = \cos \theta_h$ and define
\(
\xi_h = B_h/\sin \theta_h
\)
The matrices in (\ref{eq:harmonicintegrator}) and (\ref{eq:harmonicintegrator2}) are then
\begin{equation}\label{eq:tildemh}
\tilde{M}_h
=
\left[ \begin{matrix}\cos \theta_h & \xi_h\sin \theta_h\\ -\xi_h^{-1}\sin \theta_h & \cos \theta_h\end{matrix}\right]
\end{equation}
and
\begin{equation}\label{eq:tildemhdos}
\tilde{M}_h^n
=
\left[ \begin{matrix}\cos (n\theta_h) & \xi_h\sin (n\theta_h)\\
-\xi_h^{-1}\sin (n\theta_h) & \cos (n\theta_h)\end{matrix}\right].
\end{equation}
From \eqref{eq:tildemh} the modified Hamiltonian is found to be
\[
\tilde H_h = \frac12 \frac{\theta_h}{h}\left( \xi_hp^2+ \xi_h^{-1}q^2\right);
\]
each numerical trajectory is contained in an ellipse \(\tilde H_h= \textrm{constant}\). True solution trajectories are of course contained in circles \(p^2+q^2= \textrm{constant}\).

Comparing \eqref{eq:tildemh} with \eqref{eq:rotation} at \(t=h\), shows that order \(\nu\) is equivalent to the requirements
\[
\xi_h = 1+\mathcal{O}(h^\nu),\qquad \frac{\theta_h}{h} = 1+\mathcal{O}(h^\nu).
\]
The proximity of \(\theta_h/h\) to 1 measures the error in the angle rotated by the numerical solution at each step. Comparing \eqref{eq:tildemhdos} with \eqref{eq:rotation} at \(t=nh\), shows that, as expected, the accumulated error in the angle grows linearly with \(n\) as more and more steps are taken with a given \(h\). On the other hand, the proximity of \(\xi_h\) to 1 governs the discrepancy between the numerical ellipses and the true circular level sets of the energy. Comparing \eqref{eq:tildemhdos} with \eqref{eq:rotation} at \(t=nh\) shows that the effect of this discrepancy remains bounded as \(n\) increases. In connection with processing, errors stemming from \(\xi_h\) being different from 1 may be removed by the change of variables \(p = \xi_h^{-1/2}P\), \(q = \xi_h^{1/2}Q\) that transforms the numerical ellipses \(\xi_hp^2 +\xi_h^{-1}q^2=\) constant into circles \(P^2+Q^2=\) constant. Errors stemming from \(\theta_h/h\) being different from 1 may not be removed by changing variables/processing. For integrators of the family \eqref{eq:meth}, in general, \(\xi_h = 1+\mathcal{O}(h^2)\) and \(\theta_h/h = 1+\mathcal{O}(h^2)\). The methods that satisfy \eqref{eq:alphaminusbeta} and have enhanced conservation of energy for quadratic problems (olive double-dot-dashed curve in Fig.~\ref{fig:order}) are precisely those with \(\xi_h = 1+\mathcal{O}(h^4)\) (ellipses of low eccentricity). Similarly the methods that satisfy \eqref{eq:alphabeta} and possess effective order four (yellow solid curve in Fig.~\ref{fig:order}) are precisely those with \(\theta_h/h = 1+\mathcal{O}(h^4)\); these methods achieve order four once they are processed by a change of variables that removes at leading order the distortion in the ellipses.

\subsection{Stability}
Stability requires that the matrices in \eqref{eq:tildemhdos} remain bounded as \(n\rightarrow\infty\) with fixed \(h\). This happens if \(\tilde M_h\) has complex eigenvalues of unit modulus (which in turn happens if and only if \(|A_h|<1\)). Stability also holds if \(M_h = \pm I\), so that \(A_h=\pm 1\) and \(B_h=C_h=0\). When \(|A_h|>1\) the powers of \(\tilde M_h\) grow exponentially; when \(A_h=\pm 1\) and \(\tilde M_h\) does not coincide with \(\pm I\), \(\tilde M_h\) is a Jordan block whose powers grow linearly.

In order to investigate the stability of the integrators of the form \eqref{eq:meth}, we first see \(A_h\) in \eqref{eq:ah} as a cubic polynomial in the variable \(z=h^2\) and study the equations \(A_h= 1\) and \(A_h=-1\). Since, for fixed \(a\) and \(b\), these equations are cubic in \(z\), each of them has three roots (possibly complex and not necessarily distinct). If \(a\) or \(b\) coincide with \(0\) or \(1/2\) (\ie\ if the integrator has less than three stages) then the number of roots is smaller than three. Of course at the end of the analysis only roots with \(z\geq 0\) will play a role because \(h\) is real; however for the time being all roots are taken into account.

The analysis of the equation \(A_h(z)=1\) is easy. There is a simple root at \(z=0\) (as it corresponds to a consistent integrator) and, by removing this root, we are left with a second-degree equation. It turns out that, for non-degenerate methods, \(A_h(z)=1\) has two real, non-zero roots; these are distinct when \(a\neq b\) and coalesce for \(a=b\).

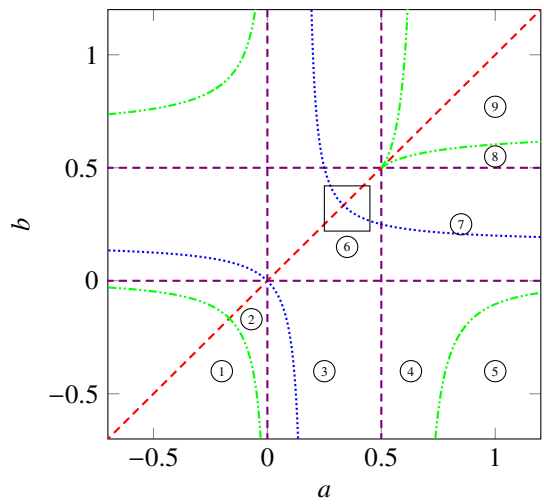
\begin{figure}
\centering
\begin{tikzpicture}
\begin{axis}[
  xmin=-0.7,
  xmax= 1.2,
  ymin=-0.7,
  ymax= 1.2,
  xlabel={\(a\)},
  ylabel={\(b\)},
  unit vector ratio*=1 1 1,
  samples=50,
]

\addplot[violet,densely dashed] coordinates {(-2.0, 0.0) (2.0,0.0)};
\addplot[violet,densely dashed] coordinates {( 0.0,-2.0) (0.0,2.0)};
\addplot[violet,densely dashed] coordinates {(-2.0, 0.5) (2.0,0.5)};
\addplot[violet,densely dashed] coordinates {( 0.5,-2.0) (0.5,2.0)};
\addplot[red,densely dashed] coordinates {(-2,-2) (2,2)};

\addplot[blue,densely dotted,domain=-1.5:0,]{degen(x)};
\addplot[blue,densely dotted,domain=-1.5:0](degen(x),x);
\addplot[blue,densely dotted,domain=0.33:2]{degen(x)};
\addplot[blue,densely dotted,domain=0.33:2](degen(x),x);

\addplot[green,densely dash dot dot,domain=-1.5:-0.2]{discrp(x)};
\addplot[green,densely dash dot dot,domain=-1.5:-0.2](discrp(x),x);
\addplot[green,densely dash dot dot,domain=0.85:2]{discrm(x)};
\addplot[green,densely dash dot dot,domain=0.85:2](discrm(x),x);
\addplot[green,densely dash dot dot,domain=0.5:2]{discrp(x)};
\addplot[green,densely dash dot dot,domain=0.5:2](discrp(x),x);
\addplot[green,densely dash dot dot,domain=-1.5:-0.15]{discrm(x)};
\addplot[green,densely dash dot dot,domain=-1.5:-0.15](discrm(x),x);

\draw (-0.20,-0.40) circle (4pt) node [scale=0.5] {1};
\draw (-0.07,-0.17) circle (4pt) node [scale=0.5] {2};
\draw (+0.25,-0.40) circle (4pt) node [scale=0.5] {3};
\draw (+0.63,-0.40) circle (4pt) node [scale=0.5] {4};
\draw (+1.00,-0.40) circle (4pt) node [scale=0.5] {5};
\draw (+0.35,+0.15) circle (4pt) node [scale=0.5] {6};
\draw (+0.85,+0.25) circle (4pt) node [scale=0.5] {7};
\draw (+1.00,+0.55) circle (4pt) node [scale=0.5] {8};
\draw (+1.00,+0.77) circle (4pt) node [scale=0.5] {9};

\coordinate (ne) at (axis cs:0.45,0.42);
\coordinate (sw) at (axis cs:0.25,0.22);
\draw[black] (sw) rectangle (ne);
\end{axis}
\end{tikzpicture}
\caption{The figure shows, for \(a\geq b\), the different possible cases that arise in the solution of the equations \(A_h(z) = \pm 1\) (see the text). The situation for \(a\leq b\) is symmetric. The framed area is described in Fig.~\ref{fig:zoom}}
\label{fig:stab}
\end{figure}

For \(A_h(z) = -1\) the situation is more involved. The discriminant of the cubic equation is
\[
ab(a+b-6ab)\Big(-12a^2b^2+8a^2b+8ab^2-6ab+\frac{1}{4}\Big)
\]
and we study the last two factors in turn. The quartic curve in the \((a,b)\) plane
\begin{equation}\label{eq:quartic}
-12a^2b^2+8a^2b+8ab^2-6ab+\frac{1}{4}=0
\end{equation}
is represented in Fig.~\ref{fig:stab} by a green double-dot-dashed line. On that line, the equation \(A_h(z) = -1\) has a double (real) root; crossing the line turns two real roots into a pair of complex conjugate pairs or vice versa. The blue dotted curved line in Fig.~\ref{fig:stab} corresponds to the hyperbola
\begin{equation}\label{eq:hyperbola}
a+b-6ab= 0.
\end{equation}
There is a double root on this line, but, at both sides of it, the double roots splits into two real roots.

We are now in a position to study the stability of the methods. Without loss of generality we consider only the set \(a\leq b\) (\(A_h\) is a symmetric function of \(a\) and \(b\)). The curve \eqref{eq:quartic} and the straight lines \(a=0\), \(a=1/2\), \(b=0\), \(b=1/2\) divide this set into nine regions, labeled in Fig.~\ref{fig:stab}. The behavior of the polynomial \(A_h(z)\) in each of these regions is shown in Fig.~\ref{fig:poly.chart}.

\begin{landscape}
\vspace*{\fill}
\begin{figure}[h!]
\centering
\begingroup
\newlength{\mylength}
\setlength{\mylength}{\textwidth}
\begin{tikzpicture}
  \begin{groupplot}[
    group style={
      group name=polychart,
      group size=3 by 3,
    },
    width=.49\mylength,
    height=.28\mylength,
    xmin=-14, ymin=-2.6,
    xmax= 34, ymax= 2.9,
    grid=both,
    grid style={line width=.2pt, draw=gray!33},
    major grid style={line width=.4pt,draw=gray!75},
    axis lines=middle,
    minor x tick num=5,
    minor y tick num=2,
    xtick={-20,-10,0,10,20,30},
    xticklabels={,,,10,20,30},
    ytick={-3,-2,-1,0,1,2,3},
    yticklabels={,,-1,,1,,},
    xlabel={\(z=h^2\)},
    ylabel={\(A_h(z)\)},
    every axis title/.append style={at={(0.5,-0.35)}},
    no markers,
    samples=40,
  ]
    \nextgroupplot[title={(1)\ \textbf{1},1,1,-1}]
      \addplot+[domain=-06:30]{A(x,-0.05,-0.80)};
    \nextgroupplot[title={(2)\ \textbf{1},-1,-1,1,1,-1}]
      \addplot+[domain=-06:34]{A(x,-0.11,-0.22)};
    \nextgroupplot[title={(3)\ 1,-1,-1,1,\textbf{1},-1}]
      \addplot+[domain=-16:04]{A(x,+0.23,-0.46)};
    \nextgroupplot[title={(4)\ -1,1,\textbf{1},-1,-1,1}]
      \addplot+[domain=-10:12]{A(x,+0.80,-0.20)};
    \nextgroupplot[title={(5)\ -1,1,\textbf{1},1}]
      \addplot+[domain=-10:12]{A(x,+0.91,-0.19)};
    \nextgroupplot[title={(6)\ \textbf{1},-1,-1,1,1,-1}]
      \addplot+[domain=-04:34]{A(x,+0.31,+0.23)};
    \nextgroupplot[title={(7)\ -1,1,\textbf{1},-1,-1,1}]
      \addplot+[domain=-10:16]{A(x,+0.80,+0.35)};
    \nextgroupplot[title={(8)\ 1,-1,-1,1,\textbf{1},-1}]
      \addplot+[domain=-10:04]{A(x,+1.10,+0.60)};
    \nextgroupplot[title={(9)\ 1,1,\textbf{1},-1}]
      \addplot+[domain=-10:04]{A(x,+1.10,+0.65)};
  \end{groupplot}
\end{tikzpicture}
\endgroup
\caption{Possible behaviors of the polynomial \(A_h(z)\). In panel number (9), the sequence 1, 1, \textbf{1}, -1 indicates that as \(z\) increases, there is first an intersection with the horizontal line \(A_h=1\), then a second intersection with \(A_h=1\), then the intersection with \(A_h=1\) at \(z=0\) (highlighted in boldface) and finally an intersection with \(A_h=-1\). The sequences in the other panels are interpreted similarly. Polynomials associated to non-degenerate methods lying on the blue dotted curve or the red dashed line of Fig.~\ref{fig:stab} are tangent to the horizontal lines \(A_h=-1\) and \(A_h=1\) (roots coalesce), respectively}
\label{fig:poly.chart}
\end{figure}
\vspace*{\fill}
\end{landscape}

%
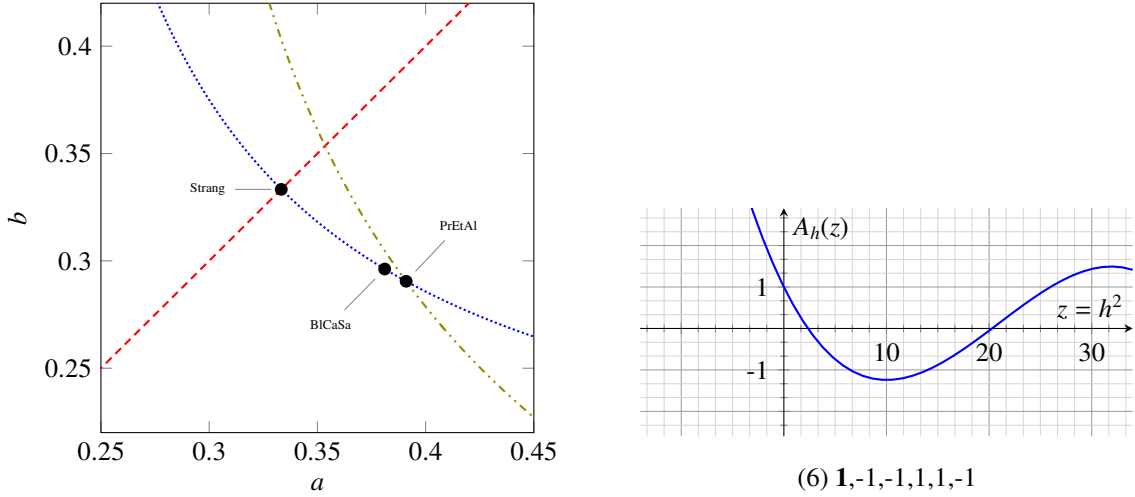
\begin{figure}
\centering
\begin{tikzpicture}
\begin{axis}[
  xmin=0.25, xmax=0.45,
  ymin=0.22, ymax=0.42,
  xlabel={\(a\)},
  ylabel={\(b\)},
  unit vector ratio*=1 1 1,
  samples=50,
]

\addplot[red,densely dashed] coordinates {(-2,-2) (2,2)};
\addplot[blue,densely dotted,domain=0.333:0.45]{degen(x)};
\addplot[blue,densely dotted,domain=0.333:0.42](degen(x),x);
\addplot[olive,dash dot dot,domain=0.268:0.42](ap(x),x);
\addplot[olive,dash dot dot,domain=0.407:0.45]{bm(x)};

\addplot[mark=*] coordinates {(0.3333, 0.3333)} node[pin=180:{\tiny Strang}]{};
\addplot[mark=*] coordinates {(0.3811, 0.2962)} node[pin=240:{\tiny BlCaSa}]{};
\addplot[mark=*] coordinates {(0.3910, 0.2905)} node[pin=60:{\tiny PrEtAl}]{};
\end{axis}
\end{tikzpicture}
\qquad
\begingroup
\setlength{\mylength}{\textwidth}
\begin{tikzpicture}
  \begin{axis}[
    width=.49\mylength,
    height=.28\mylength,
    xmin=-14, ymin=-2.6,
    xmax= 34, ymax= 2.9,
    grid=both,
    grid style={line width=.2pt, draw=gray!33},
    major grid style={line width=.4pt,draw=gray!75},
    axis lines=middle,
    minor x tick num=5,
    minor y tick num=2,
    xtick={-20,-10,0,10,20,30},
    xticklabels={,,,10,20,30},
    ytick={-3,-2,-1,0,1,2,3},
    yticklabels={,,-1,,1,,},
    xlabel={\(z=h^2\)},
    ylabel={\(A_h(z)\)},
    every axis title/.append style={at={(0.5,-0.35)}},
    no markers,
    samples=40,
    title={(6)\ \textbf{1},-1,-1,1,1,-1},
  ]
    \addplot+[domain=-04:34]{A(x,+0.31,+0.23)};
  \end{axis}
\end{tikzpicture}
\endgroup
\caption{Left panel: \((a,b)\) plane near Strang's method. The olive double-dot-dashed line corresponds to methods with high-order conservation of energy in the quadratic case. The blue dotted curve (hyperbola \eqref{eq:hyperbola}) corresponds to methods with \lq long\rq\ stability intervals (see the text). The intersection of the hyperbola and the diagonal \(a=b\) yields the longest stability interval (Strang). Right panel: behavior of the stability polynomial for methods in the region near Strang's (as in panel (6) of Fig.~\ref{fig:poly.chart}). In general, there are two intersections with the horizontal line \(A_h = -1\) near \(z=9\), \ie \(h=3\), and other two intersections with \(A_h = 1\) passed \(z=16\), \ie \(h=4\). However, polynomials associated to methods on the blue dotted curve and the red dashed line are, respectively, tangent to \(A_h = -1\) and \(A_h = 1\) near \(h=3\) and \(h=5\). This can be observed in Fig.~\ref{fig:stabilitypolynomial}}
\label{fig:zoom}
\end{figure}
Fig.~\ref{fig:zoom} focuses on an area of the \((a,b)\) plane of particular interest, as it contains methods that are not far away from Strang's splitting. From panel (6) in Fig.~\ref{fig:poly.chart} we see that, in this area, as \(z\) increases from 0, stability is typically lost at an intersection with \(A_h=-1\). However, for methods on the hyperbola \eqref{eq:hyperbola}, the graph of \(A_h(z)\) is tangent to \(A_h=-1\) and stability is, in general, lost at an intersection with \(A_h=1\).\footnote{It may be proved, see \eg\ \cite{BlCaSa}, that for the double root of \(A_h(z)=-1\), the amplification matrix is \(-I\).} The length of the stability interval is then a \emph{discontinuous} function of \(a\) and \(b\); the hyperbola \eqref{eq:hyperbola} corresponds to methods whose stability intervals are substantially longer than those of nearby methods not on the hyperbola. Here are three integrators on the hyperbola:

\begin{itemize}
\item The reference Predescu {\it et al.}  \cite{pred} suggests, in the particular case of the Newton equations, the method on the hyperbola that possesses enhanced conservation of energy in the quadratic case. The coefficients are:
\begin{equation}\label{eq:methpretal}
\textrm{PrEtAl:}\qquad a \approx 0.391 ,\qquad b \approx 0.290\:.
\end{equation}

\item At the intersection of the hyperbola with the diagonal \(a=b\), there is a double root of the equation \(A_h(z)=1\) and stability is lost, not by crossing the line \(A_h=1\), but further to the right by crossing the line \(A_h=-1\). The corresponding coefficients
\begin{equation}\label{eq:s}
\textrm{Strang:}\qquad a = 1/3 ,\qquad b = 1/3\: ,
\end{equation}
as noted before, correspond to using the Strang method \eqref{eq:vv} with step length \(h/3\). It is well known (see \eg\ \cite{BlCaSa}) that this choice of \(a\) and \(b\) leads to the longest possible stability interval \(0<h<6\) (\ie\ \(0\leq h/3 \leq 2\), where 2 is the well-known Verlet stability limit).

\item Blanes, Casas and Sanz-Serna suggested in \cite{BlCaSa} the method
\begin{equation}\label{eq:bcs}
\textrm{BlCaSa:}\qquad a \approx 0.381 ,\qquad b \approx 0.296\:.
\end{equation}
This is on the hyperbola, ensuring a favorable stability interval, and is derived by minimizing, for \(0<h<3\), the expectation of the energy error when the variables \(q\), \(p\) of the harmonic oscillator are seen as random with a Maxwell-Boltzmann distribution. This is to be compared with PrEtAl \eqref{eq:methpretal}, based on minimizing
the energy error in the limit \(h\rightarrow 0\).
\end{itemize}

The segment on the hyperbola with ends at Strang \eqref{eq:s} and PrEtAl \eqref{eq:methpretal} consists of integrators of much potential interest in molecular dynamics. When moving from PrEtAl \eqref{eq:methpretal} to Strang \eqref{eq:s}, the interval of \emph{stability improves} \cite{BlCaSa} and the size of the error constant in the \emph{energy increment} \(\Delta\) for quadratic problems \emph{worsens}. Therefore, if the user is prepared to work with smaller values of \(h\), then methods close to PrEtAl \eqref{eq:methpretal} are better, while, if \(h\) has to be taken large in order to reduce the computational effort, then methods close to Strang \eqref{eq:s} are to be preferred. It would be possible to adjust, for a given problem and a user-specified \(h\), the values of \(a\) and \(b\) along this segment of hyperbola to optimize the performance of the integrator, as it has been done in \cite{elena} in a similar situation concerning a different family of splitting integrators.

\begin{figure}
\centering
\begin{tikzpicture}
  \begin{axis}[
    width= .66\textwidth,
    height=.30\textwidth,
    xmin=-4, ymin=-2.9,
    xmax=36.5, ymax= 2.9,
    grid=both,
    grid style={line width=.2pt, draw=gray!33},
    major grid style={line width=.4pt,draw=gray!75},
    axis lines=middle,
    minor x tick num=4,
    minor y tick num=2,
    no markers,
    samples=43,
    xticklabels={,,,10,,20,,30},
    ytick={-3,-2,-1,0,1,2,3},
    yticklabels={,,-1,,1,,},
    xlabel={\(z=h^2\)},
    ylabel={\(A_h(z)\)},
    legend pos=outer north east,
  ]
    \addplot+[red,mark=none,domain=-4:38]{A(x, 0.3333, 0.3333)};
    \node[color=red] at (36,-1) {\pgfuseplotmark{*}};
    \addplot+[blue,mark=none,domain=-4:38]{A(x, 0.3811, 0.2962)};
    \addplot+[olive,mark=none,domain=-4:38]{A(x, 0.3910, 0.2905)};
    \node[color=olive] at (21.0130,1) {\pgfuseplotmark{*}};
    \addplot+[yellow,mark=none,domain=-4:38]{A(x,-0.1756,-0.1756)};
    \node[color=yellow] at (32.4330,-1) {\pgfuseplotmark{*}};
    \addplot+[orange,mark=none,domain=-4:04]{A(x,-0.1756, 1.3512)};
    \node[color=orange] at (2.4743,1) {\pgfuseplotmark{*}};
    \node[color=blue] at (21.7342,1) {\pgfuseplotmark{*}};
    \legend{Strang,BlCaSa,PrEtAl,LoSaSk,Yoshida};
  \end{axis}
\end{tikzpicture}
\caption{The polynomial \(A_h(z)\) for different methods. Note the short stability interval for the splitting of Yoshida. Methods close to PrEtAl or BlCaSa but not on the hyperbola \eqref{eq:hyperbola} lose stability by crossing \(A_h(z) = -1\) near \(z=9\) and, accordingly have stability intervals with length \(\approx 3\). For Strang, PrEtAl and BlCaSa the equation \(A_h=-1\) has a double root near \(z=9\); while PrEtAl and BlCaSa lose stability by \(A_h\) becoming larger than 1 near \(z=20\), Strang becomes unstable at \(z=36\) where \(A_h\) becomes smaller than \(-1\). LoSaSk was derived by imposing effective order 4 and that the equation \(A_h=1\) has a double root; stability is lost by \(A_h\) becoming less than \(-1\) near \(z=30\)}
\label{fig:stabilitypolynomial}
\end{figure}
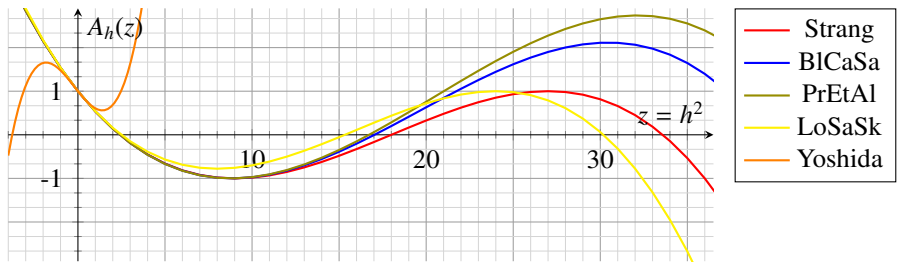

\begin{table}
\centering
\begin{tabular}{@{\quad}l|r@{.}l|r@{.}l|r@{\ \textbf{1}\ }ccr|r@{.}l@{\quad}}
\toprule
Integrator
&\multicolumn2{c|}{\(a\)}
&\multicolumn2{c|}{\(b\)}
&\multicolumn4{c|}{\(\pm1\) values}
&\multicolumn2c{\(h_{\max}\)}\\
\midrule
Strang  &  0&333333333333333 &  0&333333333333333 &     &-1&             1 &\underline{-1}& 6&000\\
BlCaSa  &  0&381119890334520 &  0&296195042611260 &     &-1&\underline{1} 1&           -1 & 4&662\\
PrEtAl  &  0&391008574596575 &  0&290485609075129 &     &-1&\underline{1} 1&           -1 & 4&584\\
LoSaSk  & -0&175603595979829 & -0&175603595979829 &     &  &             1 &\underline{-1}& 5&695\\
Yoshida & -0&175603595979829 &  1&351207191959658 & -1 1&  &               & \underline{1}& 1&573\\
\bottomrule
\end{tabular}

\caption{Coefficients and length of the stability intervals \((0,h_{\max})\) for different integrators. For the column \lq \(\pm 1\) values,\rq\ see the text} \label{tab:stability}
\end{table}

Fig.~\ref{fig:stabilitypolynomial} graphs \(A_h\) as a function of \(z\) for different methods. Table~\ref{tab:stability} provides the coefficients and the length of the stability intervals. The meaning of the column \lq\(\pm 1\) values\rq\ is as follows. In the first row, the sequence {\bf -1}, -1, 1, \underline{-1} means that as \(z\) increases, \(A_h(z)\) reaches the value 1 at \(z=0\) (consistency, boldface), then the value -1, then the value 1 and stability is lost at an intersection with \(A_h=-1\) (underline). The other rows are interpreted similarly.

For a processed method, the amplification matrix \(\widetilde{\widehat M}_h\) is given by \( P_h^{-1} \tilde M_h P_h\), where \(P_h\) is the matrix associated with \(\chi_h\) (see \eqref{eq:psihat}). Thus the eigenvalues of \(\widetilde{\widehat M}_h\) are the same as those of \(\tilde M_h\) and the processed and unprocessed integrator are simultaneously stable or unstable.

The method of effective order four with maximal stability interval corresponds to the method LoSaSk in \eqref{eq:LoSaSk} and was first found in \cite{maximal}. The following argument shows that the condition \(a=b\) leads in fact to the longest stability interval in this case. In region (1) of Fig.~\ref{fig:stab}, as \(z\) increases from 0 with \(a>b\), stability is lost when crossing the line \(A_h=1\). For methods with \(a=b\), \(A_h=1\) has a double root (the amplification matrix is the identity); the line \(A_h=1\) is touched but not crossed and stability is lost by crossing the line \(A_h=-1\) at a larger value of \(z\) (see Fig.~3), so that methods on the diagonal have stability intervals
substantially longer than methods close to the diagonal.

\section{Numerical experiments}
\label{sec:numerics}

We now present some simulations to compare the fourth-order method Yoshida \eqref{eq:yoshida}, the method LoSaSk \eqref{eq:LoSaSk} (with effective order four) and the second-order methods PrEtAl \eqref{eq:methpretal}, Strang \eqref{eq:s} and BlCaSa \eqref{eq:bcs}. We emphasize that using the three-stage scheme \eqref{eq:s} is entirely equivalent to using the standard Strang splitting \eqref{eq:vv} with step length \(h/3\). In this way we may think that we are comparing four three-stage integrators with the one-stage Strang formula while running the latter with smaller step lengths, so as to \emph{equalize the computational cost of the five algorithms}.

The integrator LoSaSk \eqref{eq:LoSaSk} is processed in \emph{all} the experiments reported to ensure results with fourth-order accuracy. As a pre-processor we use the following simple (Euler) approximation to the flow of
\eqref{eq:lambda}:
\begin{equation}\label{eq:chi1}
\chi_h(x) \approx x+h^2 \lambda\, [A,B],
\end{equation}
with the commutator evaluated at \(x\). For post-processing we take
\begin{equation}\label{eq:chi2}
\chi_h^{-1}(X) = X-h^2\lambda\, [A,B],
\end{equation}
with the commutator evaluated at \(X\). As mentioned before, \eqref{eq:lambda} may be replaced by more sophisticated choices like \eqref{eq:sophisticated} to further enhance the performance of the processed integrator; we have not considered that possibility here to avoid the need to compute the three-letter commutators needed. On the other hand, replacing the flow of \eqref{eq:lambda} and its inverse with the simple Euler approximations above implies that the processed integrator is not symplectic any longer; this should not be too worrying as the long-time properties of the implemented algorithm are really those of the unprocessed formula, which is symplectic (but more on this later).

We employ two very different test problems, one from partial differential equations and the other from molecular simulations.

\subsection{The cubic Schr\"{o}dinger equation}

Our first set of experiments concerns the well-known partial differential equation
\begin{equation}
\label{eq:nls}
i\, u_t + \partial_{xx} u + |u|^2 u = 0,
\end{equation}
in the whole real line \(-\infty<x<\infty\) or subject to periodic boundary conditions.
The equation may be written in Hamiltonian form with the kinetic and potential energy given by
\[ T(u):=\int \tfrac12|\partial_x u|^2\dx, \qquad V(u):=\int \tfrac14 |u|^4\dx \,,\]
and we split it in the two systems
\begin{equation}
\label{eq:cs:split}
u_t= i\,\partial_{xx}u, \qquad u_t = i\, |u|^2\, u,
\end{equation}
that we will identify with the symbols \(\calT\) and \(\calV\) respectively. The split system \(\calT\) may be solved in closed form in Fourier space and
the solution of \(\calV\) is clearly \(\exp(it|u_0(x)|^2)u_0(x)\), where \(u_0(x)= u(x,0)\). For processing we need the commutator; a simple computation yields
\[
[\calT,\calV] = \left[\begin{matrix}q^2+3\,p^2&-2\,qp\\-2\,qp&3\,q^2+p^2\end{matrix}\right]\,
\partial_{xx} \left[\begin{matrix}q\\p\end{matrix}\right]
- \partial_{xx} \left((q^2+p^2)\, \left[\begin{matrix}q\\p\end{matrix}\right]\right).
\]
Here, \(q\) and \(p\) are real, \(u=q+ip\) and the first and second rows of the right hand-side provide the real and imaginary part of the complex-valued commutator \([\calT,\calV]\).

We use two analytic solutions:
\begin{itemize}
\item \emph{Breather.} The expression is
\[ u =
  \left(
    \frac{b^2 \cosh(\theta) + i\, b\, \sqrt{2-b^2}\; \sinh(\theta)}
         {\cosh(\theta)-\frac{\sqrt{2-b^2}}{\sqrt{2}}\; \cos(a\, b\, x)}
    - 1
  \right)\;
  a\; \exp(i\, a^2\, t),
  \quad
  \theta=a^2\,b\,\sqrt{2-b^2}\;t, \]
where \(a\) is the limit amplitude as \(t\rightarrow \infty\) and \(b\) governs the modulation. In the experiments \(a=b=1\) and the integration is carried out in \(0\leq t\leq 3\), \(-\pi\leq x\leq \pi\) with periodic boundary conditions.
\item \emph{Soliton.} If \(a\) and \(c\) denote the amplitude and the velocity, this has the expression
\[ u = \sqrt{2a}\;\sech\left(\sqrt{a}\,(x-c\,t)\right)\;\exp\left( i\,\left(\tfrac12c\,x-(\tfrac14c^2-a)\,t\right)\right) \,. \]
Below we take \(a=2\), \(c=3\) and \(0\leq t\leq 6\). The experiments were performed in a large interval \(-20 \leq x\leq 20\) with artificially introduced periodic boundary conditions; this is possible because the soliton is exponentially small at the boundary throughout the simulation.
\end{itemize}

For numerical purposes the operator \(\partial_{xx}\) has to be discretized; two possibilities have been considered.

\begin{figure}
\centering
\begin{tikzpicture}
  \begin{groupplot}[
    group style={
      group name=pschart,
      group size=2 by 2,
    },
    group/x descriptions at=edge bottom,
    group/y descriptions at=edge left,
    xmode=log,
    ymode=log,
    xlabel={\(h\)},
    ylabel={error},
    xmin=0.005, xmax=0.06,
    cycle list={
      {red,    mark=*        },
      {blue,   mark=pentagon*},
      {olive,  mark=diamond* },
      {yellow, mark=square*  },
      {orange, mark=triangle*},
    },
    table/col sep=comma,
    table/x index={0},
  ]
    \nextgroupplot[
      ymin=0.0000012, ymax=0.1,
      legend to name={CommonLegend},
      legend style={legend columns=5},
      legend entries={Strang,BlCaSa,PrEtAl,LoSaSk,Yoshida},
    ]
      \foreach \i in {1, ..., 3, 6, 4} {
        \addplot table[y index={\i}]{brea_ps_k.csv};
      }
    \nextgroupplot[
      ymin=0.0000012, ymax=0.1,
    ]
      \foreach \i in {1, ..., 3, 6, 4} {
        \addplot table[y index={\i}]{brea_ps_d.csv};
      }
    \nextgroupplot[
      ymin=0.0000008, ymax=0.1
    ]
      \foreach \i in {1, ..., 3, 6, 4} {
        \addplot table[y index={\i}]{whit_ps_k.csv};
      }
    \nextgroupplot[
      ymin=0.0000008, ymax=0.1,
    ]
      \foreach \i in {1, ..., 3, 6, 4} {
        \addplot table[y index={\i}]{whit_ps_d.csv};
      }
  \end{groupplot}
  \node[below] at (current bounding box.south) {\ref{CommonLegend}};
\end{tikzpicture}
\caption{Error in \(u\) as a function of \(h\) for the cubic Schr\"{o}dinger equation with pseudo-spectral space discretization. The graphs in the top row correspond to the breather solution; the bottom row is for the soliton. On the left, the system \(\calT\) plays the role of system A; on the right, the system \(\calT\) plays the role of system B}
\label{fig:conv.chart.ps}
\end{figure}
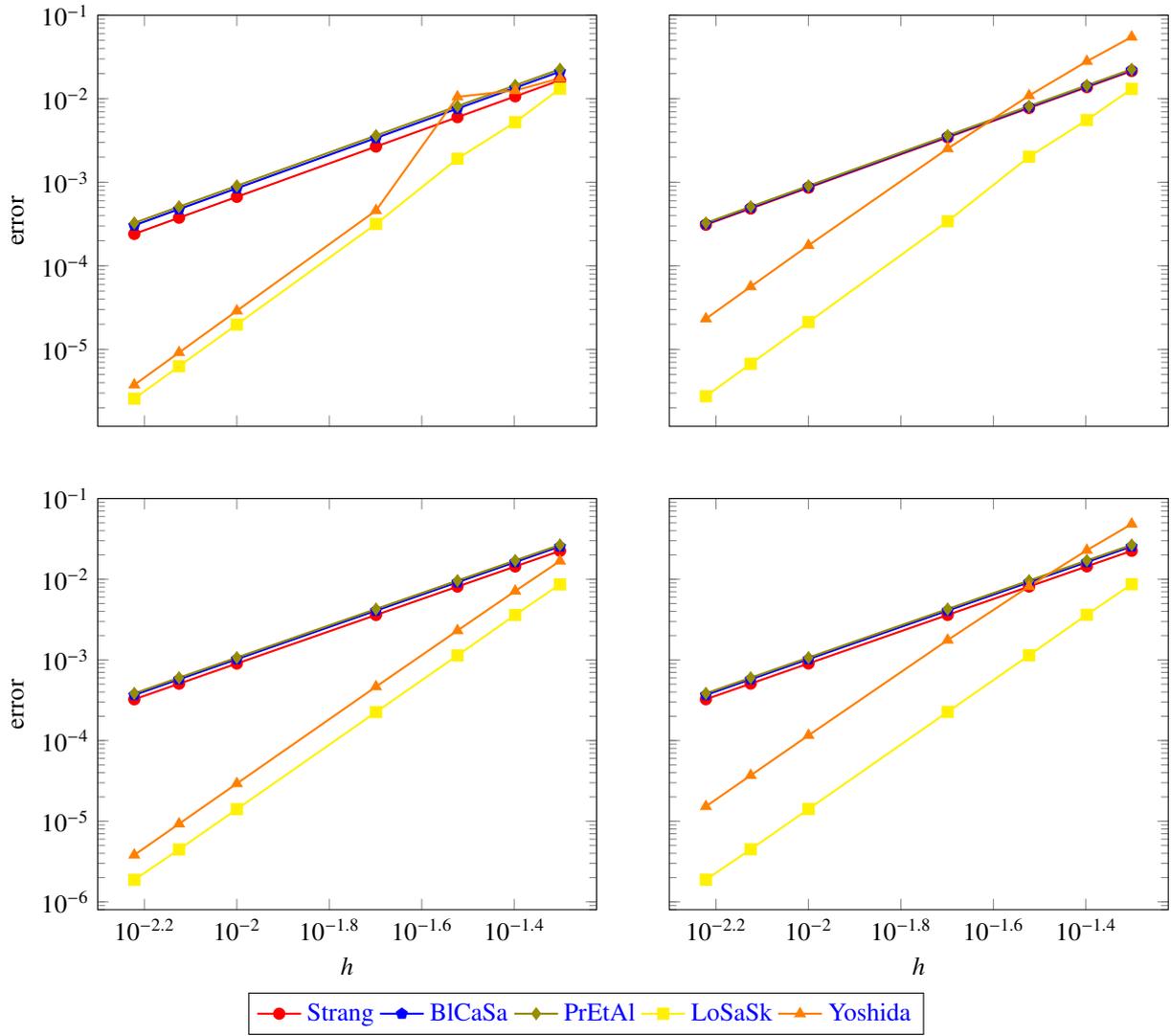

\begin{figure}
\centering
\begin{tikzpicture}
  \begin{groupplot}[
    group style={
      group name=fdchart,
      group size=2 by 2,
      xlabels at=edge bottom,
      ylabels at=edge left,
    },
    group/x descriptions at=edge bottom,
    group/y descriptions at=edge left,
    xmode=log,
    ymode=log,
    xlabel={\(h\)},
    ylabel={error},
    xmin=0.0085, xmax=0.12,
    cycle list={
      {red,    mark=*        },
      {blue,   mark=pentagon*},
      {olive,  mark=diamond* },
      {yellow, mark=square*  },
      {orange, mark=triangle*},
    },
    table/col sep=comma,
    table/x index={0},
  ]
    \nextgroupplot[
      ymin=0.0001, ymax=5,
      legend to name={CommonLegend},
      legend style={legend columns=5},
      legend entries={Strang,BlCaSa,PrEtAl,LoSaSk,Yoshida},
    ]
      \foreach \i in {1, ..., 3, 6, 4} {
        \addplot table[y index={\i}]{brea_fd_k.csv};
      }
    \nextgroupplot[
      ymin=0.0001, ymax=5,
    ]
      \foreach \i in {1, ..., 3, 6, 4} {
        \addplot table[y index={\i}]{brea_fd_d.csv};
      }
    \nextgroupplot[
      ymin=0.008,ymax=5
    ]
      \foreach \i in {1, ..., 3, 6, 4} {
        \addplot table[y index={\i}]{whit_fd_k.csv};
      }
    \nextgroupplot[
      ymin=0.008,ymax=5,
    ]
      \foreach \i in {1, ..., 3, 6, 4} {
        \addplot table[y index={\i}]{whit_fd_d.csv};
      }
  \end{groupplot}
  \node[below] at (current bounding box.south) {\ref{CommonLegend}};
\end{tikzpicture}
\caption{Error in \(u\) as a function of \(h\) for the cubic Schr\"{o}dinger equation with finite-difference space discretization. The graphs in the top row correspond to the breather solution; the bottom row is for the soliton. On the left, the system \(\calT\) plays the role of system A; on the right, the system \(\calT\) plays the role of system B}
\label{fig:conv.chart.fd}
\end{figure}
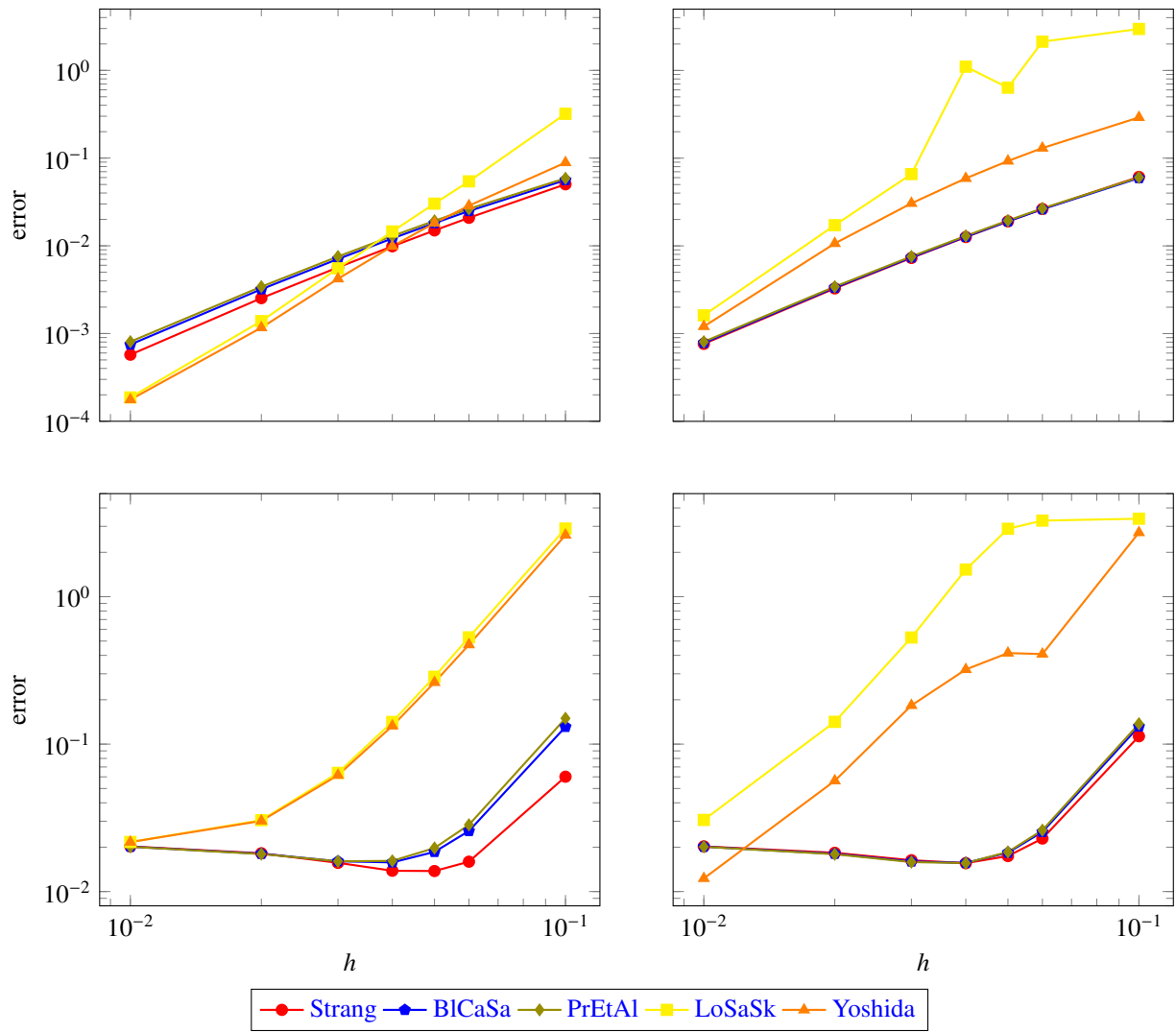

\subsubsection{Pseudo-spectral discretization}
Fig.~6 shows \(L^2\) errors with respect to the PDE solution, at the final time, when space is discretized pseudo-spectrally on a grid with 512 nodes (see \cite{fourier} for implementation details). The values of \(h\) are from the sequence 0.0500, 0.0400, 0.0300, 0.0200, 0.0100, 0.0075, 0.0060. There are four panels. The two at the top correspond to the breather and the two at the bottom to the soliton. On the left, the integrators are used with \(\calT\) and \(\calV\) in the roles of system A and system B respectively: on the right the roles are reversed.
It is well known that, for smooth solutions, pseudo-spectral techniques yield errors that are exponentially small as \(\Delta x\rightarrow 0\), \ie\ spatial errors are negligible. It is therefore of clear interest to be able to perform the time integration accurately. The four panels in the figure clearly indicate that order four is more efficient than order two.
The performance of the three second-order integrators is essentially the same. The behavior of Yoshida's formula changes markedly from left to right. For the breather, the choice of \(\calT\) as A results in errors one order of magnitude smaller than if \(\calT\) plays the role of B. For the soliton the picture is reversed. Of the five methods under consideration, the processed integrator is the best in the four panels, in spite of the crude processing applied.

\subsubsection{Finite differences}

While finite difference discretizations cannot compete with pseudo-spectral methods when integrating the cubic Schr\"{o}dinger and many other partial differential equations, it is useful to consider them here, as they are applicable to many problems that cannot be treated by spectral techniques. Fig.~7 is as Fig.~6, but now \(\partial_{xx}\) is discretized by standard central differences on a grid with 2048 nodes. For the breather this yields \(\Delta x \approx 0.003\) and for the soliton a much large value
\(\Delta x \approx 0.020\). The time-step \(h\) is taken from the sequence 0.1000, 0.0600, 0.0500, 0.0400, 0.0300, 0.0200, 0.0100.

For the soliton, the errors are always larger than \(\approx 0.02\) due to the contribution of space (the grid is not very fine and we are using second-order differences). The second order time-integrators are to be preferred; the benefits of fourth order in time would only arise for \(h\) so small that the temporal error is negligible with respect to the spatial error. The two simulations with order four perform in a similar way on the left, but in the lower right panel, Yoshida gives errors that are smaller than those of LoSaSk by a factor of three; this is of no consequence because, as we just mentioned, order four does not pay in these simulations.

For the breather, on the right panel the methods with order two are again to be preferred. On the left panel there is a cross-over situation; order two in time is advantageous for larger values of \(h\).

\subsection{A small molecule}

We now test the integrators in the simulation of the Newton equations for a molecule of the alkane
\(\rm C_{9}H_{20}\). The model under consideration has been used in \cite{cances} and \cite{extra} to compare different Markov chain samplers.
It has \(3\times 9=27\) degrees of freedom because the hydrogen atoms are lumped to the corresponding carbon atoms.
The forces are expensive to evaluate, as the model includes contributions related to bond-lengths, bond-angles, dihedral angles and also pairwise Lennard-Jones interactions. The simulations here use the parameter values in
\cite{cances} with inverse temperature \(\beta = 1\) in the system of units used there.

\begin{figure}
\centering
\begin{tikzpicture}
  \begin{groupplot}[
    group style={
      group name=pschart,
      group size=2 by 1,
    },
    group/x descriptions at=edge bottom,
    group/y descriptions at=edge left,
    xmode=log,
    ymode=log,
    xlabel={\(h\)},
    ylabel=error,
    cycle list={
      {red,    mark=*        },
      {blue,   mark=pentagon*},
      {olive,  mark=diamond* },
      {yellow, mark=square*  },
      {orange, mark=triangle*},
    },
    table/col sep=comma,
    table/x index={0},
  ]
    \nextgroupplot[
      xmin=0.0008, xmax=0.075,
      ymin=1.5e-9, ymax=0.6,
      legend entries={Strang,BlCaSa,PrEtAl,LoSaSk,Yoshida},
      legend to name={CommonLegend},
      legend style={legend columns=5},
    ]
      \foreach \i in {1, ..., 5} {
        \addplot table[y index={\i}]{alk_k.csv};
      }
    \nextgroupplot[
      xmin=0.0008, xmax=0.075,
      ymin=1.5e-9, ymax=0.6,
    ]
      \foreach \i in {1, ..., 5} {
        \addplot table[y index={\i}]{alk_d.csv};
      }
  \end{groupplot}
  \node[below] at (current bounding box.south) {\ref{CommonLegend}};
\end{tikzpicture}
\caption{Energy errors as a function of \(h\) for the alkane molecule. On the left \(T\) plays the role of system A; on the right \(T\) plays the role of system B}
\label{fig:moldyn}
\end{figure}
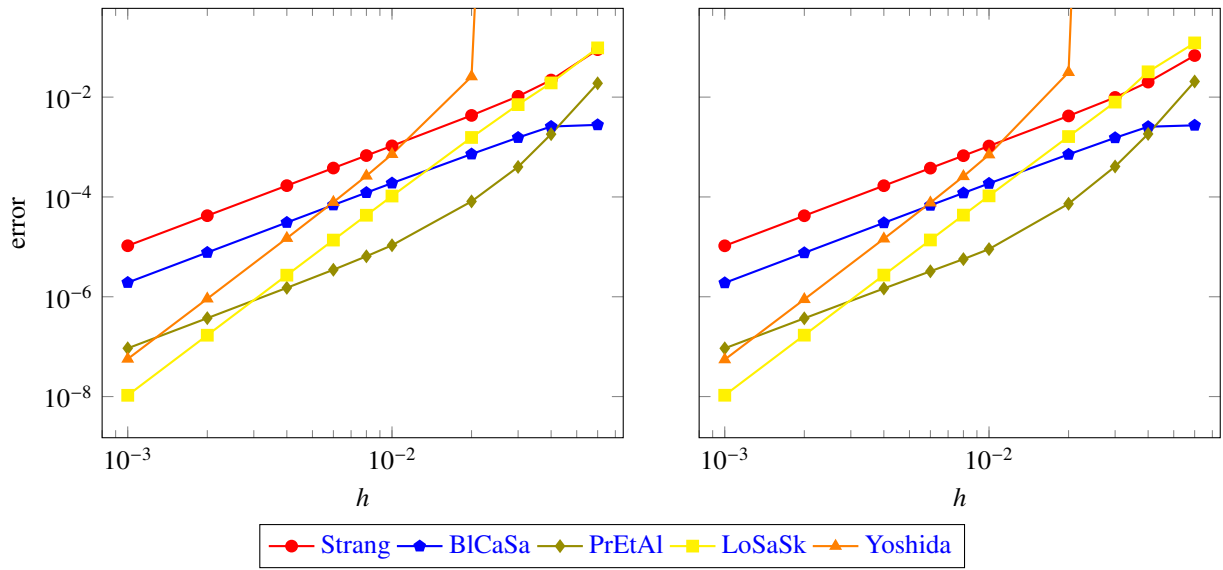

We first simulate the dynamics of the molecule in the interval \(0\leq t\leq 0.48\). The initial \(q\) corresponds to the most stable configuration of the molecule and the initial momentum \(p\) was generated randomly from the relevant Maxwell distribution. Fig.~\ref{fig:moldyn} shows the relative energy error at the final time; the results displayed correspond to averages over ten draws of the initial \(p\). On the left, \(\calT\) plays the role of system A; on the right, \(\calT\) plays the role of system B. On both panels Yoshida's method cannot operate when \(h\) is large, due to its short stability interval; as a consequence the method is of little use unless errors below \(0.1\%\) are of interest. The performance of the LoSaSk integrator is in both panels better than those of Yoshida's and Strang (which is, in the present circumstances, velocity Verlet on the left and position Verlet on the right). Both PrEtAl \eqref{eq:methpretal} and BlCaSa \eqref{eq:bcs} also improve on Strang/Verlet. Between PrEtAl \eqref{eq:methpretal} and BlCaSa \eqref{eq:bcs}, the former is more efficient for small \(h\) and the latter for larger \(h\); this is consistent with the rationale behind their derivation, as PrEtAl \eqref{eq:methpretal} is based on the behavior as \(h\rightarrow 0\) and BlCaSa \eqref{eq:bcs} on the behavior for finite \(h\).

Finally, we ran the schemes in conjunction with the Hybrid Monte Carlo or Hamiltonian Monte Carlo (HMC) method (see \eg\ \cite{cetraro}) to sample from the Boltzmann distribution corresponding to the potential energy. HMC has become very popular in Bayesian statistics and other fields of application. It performs molecular dynamics simulations to propose a state
\((q,p)\) which is accepted or rejected according to a Metropolis test based on the energy error \(\Delta\). Rejected steps are not welcome; when a rejection occurs, the Markov chain does not move and this has a bad effect on the quality of the sampling (autocorrelation increases). Acceptance is fostered by using smaller values of \(h\) leading to smaller energy errors. However, very small values of \(h\) cannot be recommended, as they would require unnecessarily expensive simulations. Theoretical results and practical experience \cite{beskos} suggest to adjust \(h\) to get acceptance rates close to \(70\%\).

Recall that processing destroys reversibility, even if the flow of \eqref{eq:lambda} could be computed exactly. In addition, the approximations
\eqref{eq:chi1}--\eqref{eq:chi2} introduce errors in the conservation of volume.
In HMC, it is important that the dynamics is simulated with an integrator that is \emph{exactly} volume-preserving and \emph{exactly} reversible\footnote{Strictly speaking these requirements may be relaxed at the expense of resorting to complicated accept/reject criteria, see \cite{fang}.} and this precludes the application of the idea of processing. Accordingly, the integrator
LoSaSk was not tested for HMC.

In the experiments the molecular dynamics legs span the interval \(0\leq t\leq 0.48\) with the initial condition chosen as above and \(\calT\) in the role of system B. We take \(h= 0.06\) which results in a satisfactory acceptance rate for the Strang splitting (which is in fact the standard position Verlet with step length \(0.02\)) and kept this value of \(h\) for all other integrators to equalize computational cost. (The value \(h=0.06\) corresponds to \(23.04\) fs, so that position Verlet runs with steps of length \(\approx 8\) fs.) The initial condition was taken as above. We generated 20 chains, each with a burn-in phase of 200 samples and a production phase of 1000 samples. Table~\ref{tab:aratio} lists the average and the standard deviation over the 20 chains of the observed acceptance rate. With the same computational cost, PrEtAl \eqref{eq:methpretal} and specially BlCaSa \eqref{eq:bcs} lead to higher acceptance rates (and therefore better sampling) than the standard Verlet integrator. Yoshida's method is not suitable for this kind of application.

\begin{table}
\centering
\begin{tabular}{@{\quad}l|r@{.}l|r@{.}l@{\quad}}
\toprule
Integrator
& \multicolumn2{c|}{\(\mu\)}
& \multicolumn2{c }{\(\sigma\)}\\
\midrule
Strang  & 77&70\% & 2&11\%\\
BlCaSa  & 96&70\% & 0&41\%\\
PrEtAl  & 91&84\% & 0&96\%\\
Yoshida &  0&00\% & 0&00\%\\
\bottomrule
\end{tabular}

\caption{Average and standard deviation of the acceptance rate in HMC sampling of the alkane Boltzmann distribution} \label{tab:aratio}
\end{table}

\section{Conclusions}
The implementation of multi-stage splitting integrators is essentially the same as the implementation of the familiar Strang/Verlet method. Therefore multi-stage formulas may be easily incorporated to software that now uses the Strang/Verlet integrator. An example is given in \cite{elena}, where multi-stage schemes have been added to an in-house version of the molecular dynamics code GROMACS.

We have studied the two-parameter family \eqref{eq:meth} of palindromic, three-stage splitting formulas.
Five specific choices of the parameters have been identified, after taking into account the order of consistency, the energy error in quadratic models and stability properties.
\begin{itemize}
\item\emph{Order four.} In order to achieve order four, the parameters in the family have to satisfy two independent equations. This yields a unique method Yoshida \eqref{eq:yoshida}, often associated with the name of Yoshida. Unfortunately the method has a very small stability interval. In addition, the application of Yoshida \eqref{eq:yoshida} is restricted to problems where the split systems \eqref{eq:A} and \eqref{eq:B} may both be integrated backward in time.

In each experiment reported in this paper, the fourth order method was outperformed by the methods below.
\item\emph{Effective order four.} There is a one-parameter subfamily of \eqref{eq:meth} consisting of methods that, while being of order two, provide approximations of order four after processing the computed solution (effective order four). Often, the processing required may be performed with very low computational cost. The application of the methods of effective order four is restricted to problems where the split systems may be integrated backward. Additionally, processing destroys the reversibility of the algorithm, and simulations with effective order four cannot be used within the HMC sampling method.

Among the one-parameter family of integrators with effective order four, we have focused on the coefficients of LoSaSk \eqref{eq:LoSaSk}, as these lead to the longest stability interval. The experiments reported suggest that, of the tested schemes, LoSaSk \eqref{eq:LoSaSk} is the best choice to integrate in time partial differential equations in situations where the spatial error is very small (\eg\ when spectral techniques are applied or finite differences/finite elements methods are used with fine spatial grids).
\item\emph{Order two.} The choice \eqref{eq:s} essentially reproduces the standard Strang/Verlet integrator with the longest stability interval. It should be used whenever one wishes to operate with values of \(h\) as large as possible, so as to save computational effort. However, if one is prepared to consider (perhaps slightly) smaller values of \(h\), it is possible to perturb the values in Strang \eqref{eq:s} to enhance the efficiency of the simulations. Perturbations have to be limited to coefficients on the hyperbola \eqref{eq:hyperbola}; leaving this curve substantially reduces the length of the stability interval. Two such perturbations, PrEtAl \eqref{eq:methpretal} and LoSaSk \eqref{eq:LoSaSk}, have been suggested in the literature, both based on the quadratic model \eqref{eq:quadratic}. The former is motivated by the wish to improve conservation of energy as \(h\rightarrow 0\) and the latter on minimizing the expectation of the energy error when the variables are random with the Maxwell-Boltzmann distribution.

The experiments reported above suggest that the second order formulas Strang \eqref{eq:s}, PrEtAl \eqref{eq:methpretal} and LoSaSk \eqref{eq:LoSaSk} are to be preferred to the integrator with effective order four in molecular dynamics problems and also in partial differential equations when the spatial discretization is not very accurate. Both PrEtAl \eqref{eq:methpretal} and LoSaSk \eqref{eq:LoSaSk} outperform the Verlet algorithm.
\end{itemize}

We have also pointed out that it would be possible, as in \cite{elena}, to adjust the coefficients of the integrator along the hyperbola \eqref{eq:hyperbola} to a given problem and a given value of \(h\), so as to use a method closer to Strang/Verlet when \(h\) is so large that stability is the main issue and to PrEtAl \eqref{eq:methpretal} for smaller values of \(h\). This point will be taken up in future research.

\section*{Acknowledgments}
JMSS has been supported by projects MTM2013-46553-C3-1-P from Ministerio de Econom\'{\i}a y Comercio, and MTM2016-77660-P (AEI/\-FEDER, UE) from Ministerio de Eco\-nom\'{\i}a, Industria y Competitividad, Spain. CMC was supported by Junta de Castilla y León, Spain, together with the European Social Fund through a postdoctoral position held at the Instituto de Matemáticas de la Universidad de Vallolid, IMUVA.

\section*{References}
\bibliographystyle{elsarticle-harv}
\bibliography{na}

\begin{thebibliography}{33}
\expandafter\ifx\csname natexlab\endcsname\relax\def\natexlab#1{#1}\fi
\expandafter\ifx\csname url\endcsname\relax
  \def\url#1{\texttt{#1}}\fi
\expandafter\ifx\csname urlprefix\endcsname\relax\def\urlprefix{URL }\fi

\bibitem[{Alamo and Sanz-Serna(2016)}]{AlSa16}
Alamo, A., Sanz-Serna, J.~M., 2016. A technique for studying strong and weak
  local errors of splitting stochastic integrators. SIAM J. Numer. Anal.
  54~(6), 3239--3257.
\newline\urlprefix\url{http://dx.doi.org/10.1137/16M1058765}

\bibitem[{Arnol$\prime$d(1997)}]{arnold}
Arnol$\prime$d, V.~I., 1997. Mathematical methods of classical mechanics.
  Vol.~60 of Graduate Texts in Mathematics. Springer-Verlag, New York,
  translated from the 1974 Russian original by K. Vogtmann and A. Weinstein,
  Corrected reprint of the second (1989) edition.

\bibitem[{Beskos et~al.(2013)Beskos, Pillai, Roberts, Sanz-Serna, and
  Stuart}]{beskos}
Beskos, A., Pillai, N., Roberts, G., Sanz-Serna, J.-M., Stuart, A., 2013.
  Optimal tuning of the hybrid {M}onte {C}arlo algorithm. Bernoulli 19~(5A),
  1501--1534.
\newline\urlprefix\url{http://dx.doi.org/10.3150/12-BEJ414}

\bibitem[{Blanes and Casas(2005)}]{negative}
Blanes, S., Casas, F., 2005. On the necessity of negative coefficients for
  operator splitting schemes of order higher than two. Appl. Numer. Math.
  54~(1), 23--37.
\newline\urlprefix\url{http://dx.doi.org/10.1016/j.apnum.2004.10.005}

\bibitem[{Blanes and Casas(2016)}]{bc}
Blanes, S., Casas, F., 2016. A concise introduction to geometric numerical
  integration. Monographs and Research Notes in Mathematics. CRC Press, Boca
  Raton.

\bibitem[{Blanes et~al.(2014)Blanes, Casas, and Sanz-Serna}]{BlCaSa}
Blanes, S., Casas, F., Sanz-Serna, J.~M., 2014. Numerical integrators for the
  hybrid {M}onte {C}arlo method. SIAM J. Sci. Comput. 36~(4), A1556--A1580.
\newline\urlprefix\url{http://dx.doi.org/10.1137/130932740}

\bibitem[{Butcher(1972)}]{Bu72}
Butcher, J.~C., 1972. An algebraic theory of integration methods. Math. Comp.
  26, 79--106.
\newline\urlprefix\url{http://dx.doi.org/10.2307/2004720}

\bibitem[{Butcher and Sanz-Serna(1996)}]{effective}
Butcher, J.~C., Sanz-Serna, J.~M., 1996. The number of conditions for a
  {R}unge-{K}utta method to have effective order {$p$}. Appl. Numer. Math.
  22~(1-3), 103--111, special issue celebrating the centenary of Runge-Kutta
  methods.
\newline\urlprefix\url{http://dx.doi.org/10.1016/S0168-9274(96)00028-1}

\bibitem[{Calvo et~al.(1994)Calvo, Murua, and Sanz-Serna}]{modified}
Calvo, M.~P., Murua, A., Sanz-Serna, J.~M., 1994. Modified equations for
  {ODE}s. In: Chaotic numerics ({G}eelong, 1993). Vol. 172 of Contemp. Math.
  Amer. Math. Soc., Providence, RI, pp. 63--74.
\newline\urlprefix\url{http://dx.doi.org/10.1090/conm/172/01798}

\bibitem[{Campos and Sanz-Serna(2015)}]{extra}
Campos, C.~M., Sanz-Serna, J.~M., 2015. Extra chance generalized hybrid {M}onte
  {C}arlo. J. Comput. Phys. 281, 365--374.
\newline\urlprefix\url{http://dx.doi.org/10.1016/j.jcp.2014.09.037}

\bibitem[{Canc\'{e}s et~al.(2007)Canc\'{e}s, Legoll, and Stoltz}]{cances}
Canc\'{e}s, E., Legoll, F., Stoltz, G., 2007. Theoretical and numerical
  comparison of some sampling methods for molecular dynamics. M2AN Math. Model.
  Numer. Anal. 41~(2), 351--389.
\newline\urlprefix\url{http://dx.doi.org/10.1051/m2an:2007014}

\bibitem[{Candy and Rozmus(1991)}]{candy}
Candy, J., Rozmus, W., 1991. A symplectic integration algorithm for separable
  {H}amiltonian functions. J. Comput. Phys. 92~(1), 230--256.
\newline\urlprefix\url{http://dx.doi.org/10.1016/0021-9991(91)90299-Z}

\bibitem[{Castella et~al.(2009)Castella, Chartier, Descombes, and
  Vilmart}]{philippe}
Castella, F., Chartier, P., Descombes, S., Vilmart, G., 2009. Splitting methods
  with complex times for parabolic equations. BIT 49~(3), 487--508.
\newline\urlprefix\url{http://dx.doi.org/10.1007/s10543-009-0235-y}

\bibitem[{Fang et~al.(2014)Fang, Sanz-Serna, and Skeel}]{fang}
Fang, Y., Sanz-Serna, J.~M., Skeel, R.~D., 2014. Compressible generalized
  hybrid monte carlo. J. Chem. Phys. 140, 174108 (10 pages).

\bibitem[{Fern\'andez-Pend\'as et~al.(2016)Fern\'andez-Pend\'as, Akhmatskaya,
  and Sanz-Serna}]{elena}
Fern\'andez-Pend\'as, M., Akhmatskaya, E., Sanz-Serna, J.~M., 2016. Adaptive
  multi-stage integrators for optimal energy conservation in molecular
  simulations. J. Comput. Phys. 327, 434--449.
\newline\urlprefix\url{http://dx.doi.org/10.1016/j.jcp.2016.09.035}

\bibitem[{Forest and Ruth(1990)}]{forest}
Forest, E., Ruth, R.~D., 1990. Fourth-order symplectic integration. Phys. D
  43~(1), 105--117.
\newline\urlprefix\url{http://dx.doi.org/10.1016/0167-2789(90)90019-L}

\bibitem[{Hairer et~al.(2010)Hairer, Lubich, and Wanner}]{hlw}
Hairer, E., Lubich, C., Wanner, G., 2010. Geometric numerical integration.
  Vol.~31 of Springer Series in Computational Mathematics. Springer,
  Heidelberg, structure-preserving algorithms for ordinary differential
  equations, Reprint of the second (2006) edition.

\bibitem[{Hansen and Ostermann(2009)}]{alexander}
Hansen, E., Ostermann, A., 2009. High order splitting methods for analytic
  semigroups exist. BIT 49~(3), 527--542.
\newline\urlprefix\url{http://dx.doi.org/10.1007/s10543-009-0236-x}

\bibitem[{Leimkuhler and Reich(2004)}]{leim}
Leimkuhler, B., Reich, S., 2004. Simulating {H}amiltonian dynamics. Vol.~14 of
  Cambridge Monographs on Applied and Computational Mathematics. Cambridge
  University Press, Cambridge.

\bibitem[{L\'opez-Marcos et~al.(1996)L\'opez-Marcos, Sanz-Serna, and
  Skeel}]{processing}
L\'opez-Marcos, M.~A., Sanz-Serna, J.~M., Skeel, R.~D., 1996. Cheap enhancement
  of symplectic integrators. In: Numerical analysis 1995 ({D}undee, 1995). Vol.
  344 of Pitman Res. Notes Math. Ser. Longman, Harlow, pp. 107--122.

\bibitem[{Lopez-Marcos et~al.(1996)Lopez-Marcos, Sanz-Serna, and
  Skeel}]{maximal}
Lopez-Marcos, M.~A., Sanz-Serna, J.~M., Skeel, R.~D., 1996. An explicit
  symplectic integrator with maximal stability interval. In: Numerical
  analysis. World Sci. Publ., River Edge, NJ, pp. 163--175.

\bibitem[{L\'opez-Marcos et~al.(1997)L\'opez-Marcos, Sanz-Serna, and
  Skeel}]{hessian}
L\'opez-Marcos, M.~A., Sanz-Serna, J.~M., Skeel, R.~D., 1997. Explicit
  symplectic integrators using {H}essian-vector products. SIAM J. Sci. Comput.
  18~(1), 223--238, dedicated to C. William Gear on the occasion of his 60th
  birthday.
\newline\urlprefix\url{http://dx.doi.org/10.1137/S1064827595288085}

\bibitem[{Murua and Sanz-Serna(1999)}]{ordercond}
Murua, A., Sanz-Serna, J.~M., 1999. Order conditions for numerical integrators
  obtained by composing simpler integrators. R. Soc. Lond. Philos. Trans. Ser.
  A Math. Phys. Eng. Sci. 357~(1754), 1079--1100.
\newline\urlprefix\url{http://dx.doi.org/10.1098/rsta.1999.0365}

\bibitem[{Predescu et~al.(2012)Predescu, Lippert, Eastwood, Ierarde, Xu,
  Jensen, Bowers, Gullingsrud, Rendleman, Dror, and Shaw}]{pred}
Predescu, C., Lippert, R.~A., Eastwood, M.~P., Ierarde, D., Xu, H., Jensen,
  M.~{\O}., Bowers, K.~J., Gullingsrud, J., Rendleman, C.~A., Dror, R.~O.,
  Shaw, D.~E., 2012. Computationally efficient molecular dynamics integrators
  with improved sampling accuracy. Mol. Phys. 110~(4), 2487 (17 pages).

\bibitem[{Sanz-Serna(1995)}]{fourier}
Sanz-Serna, J.~M., 1995. Fourier techniques in numerical methods for
  evolutionary problems. In: Computational physics ({G}ranada, 1994). Vol. 448
  of Lecture Notes in Phys. Springer, Berlin, pp. 145--200.
\newline\urlprefix\url{http://dx.doi.org/10.1007/3-540-59178-8_31}

\bibitem[{Sanz-Serna(1996)}]{fields}
Sanz-Serna, J.~M., 1996. Backward error analysis of symplectic integrators. In:
  Integration algorithms and classical mechanics ({T}oronto, {ON}, 1993).
  Vol.~10 of Fields Inst. Commun. Amer. Math. Soc., Providence, RI, pp.
  193--205.

\bibitem[{Sanz-Serna(1997)}]{geometric}
Sanz-Serna, J.~M., 1997. Geometric integration. In: The state of the art in
  numerical analysis ({Y}ork, 1996). Vol.~63 of Inst. Math. Appl. Conf. Ser.
  New Ser. Oxford Univ. Press, New York, pp. 121--143.

\bibitem[{Sanz-Serna(2014)}]{cetraro}
Sanz-Serna, J.~M., 2014. Markov chain {M}onte {C}arlo and numerical
  differential equations. In: Current challenges in stability issues for
  numerical differential equations. Vol. 2082 of Lecture Notes in Math.
  Springer, Cham, pp. 39--88.
\newline\urlprefix\url{http://dx.doi.org/10.1007/978-3-319-01300-8_2}

\bibitem[{Sanz-Serna(2016)}]{sirev}
Sanz-Serna, J.~M., 2016. Symplectic {R}unge-{K}utta schemes for adjoint
  equations, automatic differentiation, optimal control, and more. SIAM Rev.
  58~(1), 3--33.
\newline\urlprefix\url{http://dx.doi.org/10.1137/151002769}

\bibitem[{Sanz-Serna and Calvo(1994)}]{ssc}
Sanz-Serna, J.~M., Calvo, M.~P., 1994. Numerical {H}amiltonian problems. Vol.~7
  of Applied Mathematics and Mathematical Computation. Chapman \& Hall, London.

\bibitem[{Schlick(2002)}]{schlick}
Schlick, T., 2002. Molecular modeling and simulation. Vol.~21 of
  Interdisciplinary Applied Mathematics. Springer-Verlag, New York, an
  interdisciplinary guide.
\newline\urlprefix\url{http://dx.doi.org/10.1007/978-0-387-22464-0}

\bibitem[{Strang(1963)}]{strang}
Strang, G., 1963. Accurate partial difference methods. {I}. {L}inear {C}auchy
  problems. Arch. Rational Mech. Anal. 12, 392--402.
\newline\urlprefix\url{http://dx.doi.org/10.1007/BF00281235}

\bibitem[{Yoshida(1990)}]{yoshida}
Yoshida, H., 1990. Construction of higher order symplectic integrators. Phys.
  Lett. A 150~(5-7), 262--268.
\newline\urlprefix\url{http://dx.doi.org/10.1016/0375-9601(90)90092-3}

\end{thebibliography}
\end{document}